\documentclass[graybox]{svmult}

\usepackage[utf8]{inputenc}
\usepackage{microtype}

\usepackage{mathptmx}       
\usepackage{helvet}         
\usepackage{courier}        
\usepackage{type1cm}        

\usepackage{amsmath,amsfonts,amssymb,bm}
\usepackage{mathrsfs}
\usepackage{bm}
\usepackage{mathtools}
\usepackage{booktabs}
\usepackage[mathscr]{eucal}

\usepackage[ruled,linesnumbered]{algorithm2e}

\usepackage{graphics}

\usepackage{caption}
\usepackage{subcaption}

\usepackage{pgfplots}
\pgfplotsset{compat=1.16}
\usepackage{pgf,tikz}

\usetikzlibrary{calc}
\usetikzlibrary{shapes.geometric,arrows,shapes,chains,automata}
\usetikzlibrary{external}
\usepackage{amsfonts}
\usepackage{wrapfig}
\usepackage[pdfborder={0 0 0}]{hyperref}

\usepackage{enumitem}

\usepackage{makeidx}

\usepackage{eurosym}

\usepackage{upquote}  

\newcommand{\R}{\mathbb R}

\usepackage{dsfont}			

\def\c{{\mathcal{C}}}
\def\e{{\mathcal{E}}}
\def\g{{\mathcal{G}}}
\def\n{{\mathcal{N}}}
\def\D{{\mathcal{D}}}

\def\oN{{\overline{N}}}
\def\on{{\overline{n}}}

\def\eins{\mathds{1}}

\begin{document}

\title*{Modelling and Simulation of District Heating Networks}
\author{C. J\"akle, L. Reichle, and S. Volkwein}

\institute{Christian J\"akle, Lena Reichle, Stefan Volkwein \at Department of Mathematics and Statistics, University of Konstanz, Konstanz, Germany. \email{{christian.jaekle,lena.reichle,stefan.volkwein}@uni-konstanz.de}}

\maketitle

\abstract{In the present paper a detailed mathematical model is derived for district heating networks. After semidiscretization of the convective heat equation and introducing coupling conditions at the nodes of the network one gets a high-dimensional system of differential-algebraic equations (DAEs). Neglecting temporal changes of the water velocity in the pipes, the numerical solutions do not change significantly and the DAEs have index one. Numerical experiments illustrate that the model describes the real situation very well.}


\section{Introduction}

The efficient use of energy, particularly renewable energy sources, plays an important role in today's discussions. Therefore district heating becomes more important in energy use, as it is flexible in the supply of different forms of energy.\\
District heating is a system that transports heat energy via a network of pipelines from a central power plant to different consumers with different requests. For a long time, district heating has been considered as a static problem, due to the volatility and diversity in energies and  the different requests of the consumers, this assumption becomes obsolete. Therefore, it has become more important to simulate the dynamic processes based on changes in supply/demand of the consumers and energy productions. This simulation need to be used over large time horizons.
The corresponding mathematical task is challenging. On the one hand, it is important to obtain an accurate modeling that can realistically simulate, on the other hand, the modeling should be efficient and implementable for the computer. Therefore, the formulation of a numerically efficient and stable model is an important building task.\\
In order to formulate this model, first of all, a suitable mathematical model of the network is needed, secondly, a system of one-dimensional nonlinear hyperbolic partial differential equations (PDEs) is needed to model the temperature and velocity flow in a pipe over time, and lastly, algebraic equations are needed, which ensure the mass conservation, pressure continuity, mixing temperature at the nodes and guarantee that the consumers requests are fulfilled. 
In this work we use a spatial discretization for the PDEs to get only ordinary differential equations (ODEs) for the temperature and velocity. With this, we get so-called differential algebraic 
equations (DAEs), which are mathematically challenging. 
DAEs are similar in some aspects to ODEs, but differ in some aspects, which make it more difficult to solve them. \\
An important role in the theory is the existence of consistent initial values and also of (unique) solutions. In this context, an index concept has been introduced. The higher the index is, the more complex the problem becomes, so it is obvious to try to keep it as low as possible.
In this paper we have made some simplifications in the differential equations in order to get a problem with a lower index, which nevertheless shows a certain accuracy compared to the original problem.\\
The article is organized as follows. In Section \ref{sec:DAE}, we give a short introduction into differential algebraic equations, in particular semi-explicit DAEs and review results for the unique existence for index 1 and index 2 problems. After that we formulate the mathematical model of the district heating network in Section \ref{sec:3}. In Section \ref{sec:4}, numerical results concerning the two different index concepts are compared and a simulation of a real network is presented. Finally, we draw some conclusions in Section \ref{sec:5}.

\section{Differential Algebraic Equation}\label{sec:DAE}
\label{sec:2}
\subsection{Preliminary Notes}\label{sec:intro}
In this section, we briefly recall the definition of a DAE and the index concept we use in this work. Furthermore, we give some solution results for semi-explicit DAEs. For further details we refer to \cite{BuchvomSchropp,gerdts2011optimal,kunkelBook,mehrmann2012index, schwarz2018new, Schwarz2000Consistent}, for instance.

\begin{definition}
	A DAE is an \emph{implicit ODE} of the form
	\begin{align}\label{generalDAE}
		F\big(t,z(t),z'(t)\big) = 0
	\end{align}
	with $F: I\times\mathcal{D}_F \rightarrow \mathbb{R}^{n}$ and $I \subset \mathbb{R}$ an open interval, $\mathcal{D}_F \subset \mathbb{R}^{n_z} \times \mathbb{R}^{n_z} $, $n_z \in \mathbb{N}$. 
\end{definition}

Let $z$ be a solution to \eqref{generalDAE} and $F$ continuously differentiable. We suppose that there exists an $\varepsilon>0$ such that the Jacobian matrix
\begin{align*}
F_{z'}(s,v,w)=\frac{\partial F}{\partial z'}(s,v,w)\in\mathbb R^{n_z\times n_z}
\end{align*}
is regular for all $(s,v,w)$ belonging to the neighborhood
\begin{align*}
\mathcal U=\Big\{(\tilde s,\tilde v,\tilde w)\in I\times\mathcal{D}_F\,\Big|\,\exists t\in I:|\tilde s-t|+{\|\tilde v-z(t)\|}_2+{\|\tilde w-z'(t)\|}_2<\varepsilon\Big\}.
\end{align*}
Then we can use the implicit function theorem and rewrite the system as a classical explicit ODE system.

From now on, we assume that $F$ is differentiable (on an open set containing $I\times\mathcal{D}_F$). Moreover, the Jacobian matrix $D_{z'}F$ is continuous on $I \times \mathcal{D}_F$ and singular for at least one point $(t, z(t), z'(t)) \in I\times\mathcal{D}_F$.

\begin{remark}
	In many cases the abstract DAE \eqref{generalDAE} has the structure
	\begin{subequations}\label{eq:semiDAE}
	\begin{align}
			x'(t) &= f(t, x(t), y(t)),\label{eq:semiDAE1}\\
			0 &= g(t,x(t),y(t))\label{eq:semiDAE2}
	\end{align}
	\end{subequations}
	with $x(t) \in \mathbb{R}^{n_x}, \, y(t) \in \mathbb{R}^{n_y}$ and functions $f: I\times \mathcal{D}_x \times
	\mathcal{D}_y \rightarrow \mathbb{R}^{n_x}$, $g: I\times \mathcal{D}_x \times \mathcal{D}_y \rightarrow
	\mathbb{R}^{n_y}$ with 
	$I \subset \mathbb{R}$ an open interval and $\mathcal{D}_x \times \mathcal{D}_y \subset \mathbb{R}^{n_x}
	\times \mathbb{R}^{n_y}$. 
	Then, we call \eqref{eq:semiDAE} a semi-explicit DAE. In the case of semi-explicit DAEs we assume from now on, that $f$ is continuous and $g$ is continuously differentiable 
	(on an open set containing $I \times \mathcal{D}_x \times \mathcal{D}_y$).\\
	In general DAEs can be classified through different index concepts. One is the differentiation index, introduced in the following. 
	The different index shows different properties, like solvability results of a DAE. \hfill$\Diamond$
\end{remark}
\begin{definition}
	Let the DAE $F(t,z(t),z'(t))=0$, $t\in I$, have a locally unique solution and let the function $F$ be 
	sufficiently often continuously differentiable in a neighborhood of the solution. To a given $m \in \mathbb{N}$ and $t\in I$ consider the equations
	\begin{align}\label{indexEqu}
		F(t,z(t),z'(t))=0, \quad \frac{\mathrm d F}{\mathrm dt}(t,z(t),z'(t))=0,\ldots,\quad\frac{\mathrm d^m F}{\mathrm dt^m}(t,z(t),z'(t)) = 0.
	\end{align}
	The smallest natural number $m$, for that \eqref{indexEqu} can be written as 
	\begin{align}\label{indexEqu2}
		z'(t) = \phi (t, z(t))\quad\text{for }t\in I
	\end{align}
	is called \emph{differentiation index} $di = m$. Equation \eqref{indexEqu2} is called the \emph{underlying ODE of the DAE} $F(t,z(t),z'(t))=0$.
\end{definition}
Note, that an explicit ODE has differentiation index $di = 0$, an algebraic equation $F(t,z(t)) = 0$ with regular Jacobian matrix $F_z(t,z(t))$ has differentiation index $di = 1$.\\
In the case of semi-explicit DAEs of the form \eqref{eq:semiDAE}, the differentiation index depends on the function $g$. 
Assume that we have a solution $(x,y)$ of a semi-explicit DAE-system \eqref{eq:semiDAE}. 
We already have a differential equation for $x$, therefore we need to calculate a differential equation for $y$. 
In this case, we differentiate $g$ one time with respect to $t$ and get
	\begin{equation}\label{calcIndex}
	\begin{aligned}
		0 &= g_t(t,x(t),y(t)) + g_x(t,x(t),y(t)) x'(t) + g_y(t,x(t),y(t)) y'(t)\\
		&= g_t(t,x(t),y(t)) + g_x(t,x(t),y(t)) f(t,x(t),y(t)) + g_y(t,x(t),y(t)) y'(t).
	\end{aligned}
	\end{equation}
If the matrix
\begin{align}\label{indexonecond}
	g_y(t,x(t),y(t))\in\mathbb R^{n_y\times n_y} \text{ is regular}
\end{align}
in a neighborhood of $(t,x(t),y(t))$ for $t\in I$, \eqref{calcIndex} is solvable for $y'$ and we get the underlying ODE of the DAE
	\begin{align*}
		x'(t)&= f(t,x(t), y(t)) \\
		y'(t)&= -g_y(t,x(t),y(t))^{-1} \big(g_t(t,x(t),y(t)) + g_x(t,x(t),y(t))f(t,x(t),y(t))\big).
	\end{align*}
Therefore the semi-explicit DAE has differential index $di = 1$. We call \eqref{indexonecond} the \emph{index one-condition}. 

DAEs with a differential index $di = 0$ or $di = 1$ are from a numerical point of view much easier to handle than DAEs with 
higher differentiation index $di \geq 2$. Therefore it is common to reduce the index if this is possible. Assume that we have a semi-explicit DAE in the 
form of \eqref{eq:semiDAE}. An easy way to reduce the index, is to use \eqref{calcIndex} instead of \eqref{eq:semiDAE2}. With this method 
one loose information and the solution of the reduced system generally does not correspond with the solution of the original DAE, therefore one need additional conditions, 
see Lemma \ref{lem:equiv}. 

Another interesting index concept is the perturbation index. The perturbation index indicates the influence of perturbations and their derivatives
on the solution and therefore addresses the stability of DAEs. For more details, we refer the reader to  \cite[Section $1.1.1$]{gerdts2011optimal} and \cite[Chapter 7]{Hairer2}.
In many cases the differentiation index correspond with the other index concepts. 

One big problem with higher index DAEs is to get consistent initial values. Compared to ODEs not every intitial value is consistent. 
The following definition is based on \cite[Section $5.3.4$]{Petzold} and \cite[Section $3.1$]{burger2017survey}.
\begin{definition}\label{def:consistent_initial_value}
For a general DAE \eqref{generalDAE} with differentiation index $d$ and for a sufficiently often continuously differentiable function $F$, the initial value $z_0 = z(t_0)$ is said to be \emph{consistent at $t_0$}, if the equation
\begin{align}
F^{(j)}\big(t_0,z_0^0,z_0^1, \ldots, z^{j+1}_0\big) = 0\quad\text{for }j = 0,\ldots, d-1 			\label{consValueEquation}
\end{align}
has a solution $(z_0^0,z_0^1, \ldots, z^{j+1}_0)\in\mathbb R^{n_z \times (j+2)}$, where $F^{(0)}(t,z(t),z'(t)) := F(t,z(t),z'(t))$ is set and
\begin{align*}
F^{(j)}(t,z^0,\ldots,z^{j+1}) :=& \frac{\partial F^{(j-1)}}{\partial t} (t,z^0,z^1,\ldots, z^j) + \sum_{l=0}^{j} \frac{\partial F^{(j-1)}}{\partial z^{(l)}} (t,z^0,\ldots, z^j)z^{l+1}
\end{align*}
holds for $(t,z^0,\ldots,z^{j+1})\in I\times \mathbb R^{n_z}\times\ldots\times\mathbb R^{n_z}$ and $j = 0,\ldots, d-1$.
\end{definition}

Note that the system of nonlinear equations \eqref{consValueEquation} in general has many solutions and additional conditions are required to obtain a 
particular consistent initial value, which might be relevant for a particular application.

Again, in case of semi-explicit DAEs, it depends on the function $g$ if an initial value is consistent or not. Assume that we have a solution $(x,y)$ of a semi-explicit DAE of the form \eqref{eq:semiDAE}.\\
In case of differentiation index $di = 1$, the initial value $(x_0, y_0) \in \mathcal{D}_x\times \mathcal{D}_y$ is said to be \emph{consistent at $t_0 \in I$}, if 
\begin{align*}
	g(t_0,x_0,y_0) = 0 \mbox{ holds.}
\end{align*}
Therefore, we define for given $t\in I$
\begin{align*}
	\mathcal{M}_0(t) := \big\{(x,y) \in \mathcal{D}_x\times \mathcal{D}_y \,\big|\, g(t,x,y) = 0\big\}
\end{align*}
the set of all consistent initial values for semi-explicit DAEs with differentiation index $di = 1$ at the starting time $t$.

In case of differentiation index $di = 2$, the initial value $(x_0, y_0) \in \mathcal{D}_x\times \mathcal{D}_y$ is said to be \emph{consistent at $t_0 \in I$}, if 
\begin{align*}
	g(t_0,x_0,y_0) = 0 
\end{align*}
holds and additionally 
\begin{align*}
	 \partial_tg(t_0,x_0,y_0) + \partial_xg(t_0,x_0,y_0)f(t_0,x_0,y_0)+\partial_yg(t_0,x_0,y_0)w_0 = 0,
\end{align*}
introduced by \eqref{calcIndex}, has a solution $w_0 \in \R^{n_y}$.
Again, we define the set of all consistent initial values 
\begin{align*}
		\mathcal{M}_1(t) := \big\{(x,y) \in \mathcal{D}_x \times \mathcal{D}_y \,\big|\,& 
		g(t,x,y) = 0, \, \exists\, w \in \mathbb{R}^{n_y}\mbox{ with}\\
		&\partial_tg(t,x,y) + \partial_xg(t,x,y)f(t,x,y)+\partial_yg(t,x,y)w = 0\big\}
\end{align*}
for $t \in I$ for semi-explicit DAEs with differentiation index $di = 2$ at the time $t$.

\subsection{Solvability Results}

In the following we only work with semi-explicit DAEs and give some solvability results for index-1 and index-2 semi-explicit DAEs. 
The following Section is oriented on \cite{Schwarz2000Consistent}. For this section we make use of the following hypothesis.

\smallskip
\noindent
\textbf{Assumption (A1)} \emph{The function $f$ is a continuous function and at least uniformly Lipschitz-continuous in $(x,y)$ and $g$ 
is differentiable with a uniformly Lipschitz-continuous 
 	derivative on an open subset $I \times \mathcal{D}_1 \times \mathcal{D}_2 \subset I \times \mathcal{D}_x \times \mathcal{D}_y$.}

\begin{theorem}
	Let the semi-explicit DAE \eqref{eq:semiDAE} have differentiation index $di = 1$ on $I \times \mathcal{D}_1 \times \mathcal{D}_2$. 
	Then for $t_0\in I$ and $(x_0,y_0) \in \mathcal{M}_0(t_0)\cap( \D_1\times\D_2)$ there exists a locally unique solution 
	$x : \tilde{I} \to \mathbb{R}^{n_x}$, $y : \tilde{I} \to \mathbb{R}^{n_y}$ in $C^1$ with $t_0 \in \tilde{I} \subset I$ an open interval and 
	$x(t_0)  = x_0$ and $y(t_0) = y_0$.
\end{theorem}

\begin{proof}
	Assume that we have $(x_0, y_0) \in \mathcal{M}_0(t_0) \cap \D_1 \times \D_2$ for $t_0 \in I$ given, then we already have consistent initial values. 
	Due to the fact that the semi-explicit DAE has differentiation index 
	$di = 1$, the partial derivative $g_y(t_0,x_0,y_0)$ has to be regular. 
	The implicit function theorem implies, that there exists an open set $U \subset I \times \mathbb{R}^{n_x}$ containing $(t_0, x_0)$ and 
	an unique coninuously differentiable function $y:U \to \R^{n_y}$ with $y(t_0,x_0) = y_0$ and $g(t,x,y(t,x)) = 0$ for all $(t,x) \in U$. 
	The Picard-Lindelöf theorem implies that there 
	exists $\tilde{I} \subset I $ an open interval with $t_0 \in \tilde{I}$ and a local unique solution $x$ of the initial value problem
	\begin{align*}
		x'(t) = f(t,x(t),y(t,x(t)))\text{ for }t\in \tilde I,\quad x(t_0)= x_0.
	\end{align*}
	With this, $(x,y)$ solves the semi-explicit DAE for the initial value $(x(t_0,y(t_0)) = (x_0,y_0)$.\hfill$\Box$
\end{proof}
\begin{lemma}\label{lem:equiv}
	Suppose that a semi-explicit DAE of the form \eqref{eq:semiDAE} has differentiation index $di = 1$ or $di = 2$ on 
	$I \times \mathcal{D}_1 \times \mathcal{D}_2 $.  
	Then the reduced problem
	\begin{subequations}\label{eq:reducesem}
		\begin{align}
			x'(t) &= f(t, x(t), y(t))\label{eq:reducesem1} \\
			0 &= \partial_tg(t,x(t),y(t)) + \partial_xg(t,x(t),y(t))f(t,x(t),y(t))\label{eq:reducesem2}\\
			\nonumber
			&\quad+\partial_yg(t,x(t),y(t))y^{\prime}(t)
	\end{align}
	\end{subequations}
	has differentiation index $di_{red} = di -1$.\\
Assume that there are given consistent initial values $(t_0, x_0, y_0) \in I \times \mathcal{D}_1 \times \mathcal{D}_2$ 
for the reduced problem \eqref{eq:reducesem} fulfilling 
	\begin{align}\label{condred}
		g(t_0,x(t_0),y(t_0)) = 0.
	\end{align}
	Then the solution of the reduced problem \eqref{eq:reducesem} for these initial values and the solution of the original 
	problem of the form \eqref{eq:semiDAE} 
	 are the same. 
\end{lemma}
\begin{proof}
	It follows directly from the definition of the differentiation index that the reduced DAE \eqref{eq:reducesem} has differentiation index $di_{red} = di-1$. To show the equivalence of the solutions, first, 
	every solution of the original DAE remains a solution of the reduced DAE. Conversely, for $\alpha(t):= g(t,x(t),y(t))$ we have 
	\begin{align*}
		\alpha'(t) &= \partial_tg(t,x(t),y(t)) + \partial_xg(t,x(t),y(t))f(t,x(t),y(t))+\partial_yg(t,x(t),y(t))y^{\prime}(t)\\
		& = 0
	\end{align*}
	and \eqref{condred} implies that $\alpha(t_0) = 0$. Therefore, $\alpha$ vanishes identically and the result follows.\hfill$\Box$
\end{proof}
The following Theorem is motivated by \cite[Theorem 1.3.3]{Schwarz2000Consistent}
\begin{theorem}
	Suppose that a semi-explicit DAE as \eqref{eq:semiDAE} has differentiation index $di = 2$ on $I \times \mathcal{D}_1 \times \mathcal{D}_2$. 
	Moreover, it holds
	\begin{align}\label{ass:ker_cons}
	\tag{\textbf{A2}}
	\ker \partial_yg(t,x,y) \mbox{ does not depend on } (t,x,y)\in I \times \mathcal{D}_1 \times \mathcal{D}_2.
	\end{align}	 
	Then for $(x_0,y_0)\in \mathcal{M}_1(t_0) \cap \D_1\times \D_2$ with $t_0 \in I$ 
	there exists a locally unique solution 
	$x : \tilde{I} \to \mathcal{D}_1$, $y : \tilde{I} \to \mathcal{D}_2$ in $C^1$ with $t_0 \in \tilde{I} \subset I$ 
	an open interval and $(x(t_0),y(t_0)) = (x_0,y_0)$.
\end{theorem}
\begin{proof}
	Assume that $(x_0, y_0) \in \mathcal{M}_1(t_0)\cap \D_1\times \D_2$ for $t_0 \in I$ is given. 
	Then we already have consistent initial values.  
	The definition of $\mathcal{M}_1(t_0)$ implies 
	that we can choose $w_0 \in \R^{n_y}$ with 
	\begin{align*}
		\partial_tg(t_0,x_0,y_0) + \partial_xg(t_0,x_0,y_0)f(t_0,x_0,y_0)+\partial_yg(t_0,x_0,y_0)w_0 = 0.
	\end{align*}
	Due to Lemma \ref{lem:equiv} and the definition of $\mathcal{M}_1(t_0)$ we can use the reduced problem of the form 
	\eqref{eq:reducesem} instead of the original problem
	 to proof the assertion. We consider
	 \begin{subequations}\label{eq:reducesem_b}
		\begin{align}
			x'(t) &= f(t, x(t), y(t))\label{eq:reducesem_b1} \\
			0 &= \partial_tg(t,x(t),y(t)) + \partial_xg(t,x(t),y(t))f(t,x(t),y(t))\label{eq:reducesem_b2}\\
			\nonumber
			&\quad+\partial_yg(t,x(t),y(t))y^{\prime}(t),\\
			x(t_0) &= x_0,\quad y(t_0) = y_0 \label{eq:reducesem_b3}
	\end{align}
	\end{subequations}
	with $t \in I$. 
	Lemma \ref{lem:equiv} implies that the reduced problem \eqref{eq:reducesem_b} has differentiation index $di_{red} = 1$. 
	Therefore $\partial_yg(t,x,y)$ is singular. 
	Due to \eqref{ass:ker_cons} we can define the projector $Q$ onto $\ker \partial_yg(t,x,y)$ and 
	$P := \eins - Q$. Then for the projections we have the following equalities:
	\begin{align}
		Q^2 &= Q \label{eq:eq_first}\\
		QP &= 0 \label{eq:eq_second}\\
		\partial_yg(t,x,y)P &= \partial_yg(t,x,y) \label{eq:eq_third}\\
		\partial_yg(t,x,y)Q &= 0 \label{eq:eq_fourth}
	\end{align}
	for all $(t,x,y) \in I \times \D_1 \times \D_2$. 
	Consider $y \in \mathcal{D}_2, \, w \in \R^{n_y}$ and define the variables 
	\begin{align*}
		u &= Py \\
		v &:= Pw+Qy 
	\end{align*}
	we call $u$ the regular variable and $v$ the singular variable. We set $u_0 := Py_0$ and $v_0 := Pw_0 + Qy_0$. 
	With \eqref{eq:eq_first} and \eqref{eq:eq_second} it follows 
	\begin{align*}
		y = Py + Qy = Py + QPw + QQy = u + Qv. 
	\end{align*}
	Due to \eqref{eq:eq_third} and \eqref{eq:eq_fourth} we can rewrite \eqref{eq:reducesem_b2}
	\begin{align*}
		&\partial_tg(t,x,y) + \partial_xg(t,x,y)f(t,x,y)+\partial_yg(t,x,y)w \\
			&= h(t,x,y) + \partial_yg(t,x,y)w \\
			&= h(t,x,u+Qv) + \partial_yg(t,x,u+Qv)Pw \\
			&=  h(t,x,u+Qv) + \partial_yg(t,x,u+Qv)Pw + \partial_yg(t,x,u+Qv)Qy \\
			&= h(t,x,u+Qv) + \partial_yg(t,x,u+Qv)v
	\end{align*}
	with the function $h : I \times \D_1 \times \D_2 \to \R^{n_y}$ given through 
	\begin{align*}
	h(t,x,y) &:= \partial_tg(t,x,y) + \partial_xg(t,x,y)f(t,x,y).
	\end{align*}
	We define the reduced function $\tilde{g}: I \times \D_1 \times \D_2 \times \R^{n_y} \to \R^{n_y}$ through
	\begin{align*}
	 \tilde{g}(t,x,u,v) &:= h(t,x,u+Qv) + \partial_yg(t,x,u+Qv)v.
	\end{align*}
	For the partial derivative of $\tilde{g}$ it holds
	\begin{align*}
		&\partial_{v}\tilde{g}(t,x,u,v)\\
		&\quad=\partial_{v}h(t,x,u+Qv) + \partial_{v}(\partial_yg(t,x,u+Qv))v + \partial_yg(t,x,u+Qv)  \\
		&\quad=\partial_y\left(\partial_tg(t,x,u+Qv)+\partial_xg(t,x,u+Qv)f(t,x,u+Qv)\right)Q \\
			&\qquad +\partial_y^2g(t,x,u+Qv)Qv + \partial_yg(t,x,u+Qv)
	\end{align*}
	Because the original DAE has differentiation index $di = 2$ it follows that 
	for $u_0,\, v_0\in \R^{n_y}$ chosen as before we have
	\begin{align*}
		 \partial_v\tilde{g}(t_0,x_0,u_0,v_0) \mbox{ is regular.}
	\end{align*}
	In particular with \eqref{eq:eq_first}, \eqref{eq:eq_second} and the definition of $u_0,\,v_0$, it follows
	\begin{align*}
		&\tilde{g}(t_0,x_0,u_0,v_0)\\
		&\quad= h(t_0,x_0,Py_0+Q(Pw_0+Qy_0)) + \partial_yg(t_0,x_0,Py_0+Q(Pw_0+Qy_0))(Pw_0+Qy_0)\\
		 &\quad=  h(t_0,x_0,Py_0+Qy_0) + \partial_yg(t_0,x_0,Py_0+Qy_0)w_0\\
		 &\quad= \partial_tg(t_0,x_0,y_0)+\partial_xg(t_0,x_0,y_0)f(t_0,x_0,y_0)+\partial_yg(t_0,x_0,y_0)w_0 = 0.
	\end{align*}
	Due to the implicit function theorem there exists a neighbourhood $U$ of $(t_0,x_0,u_0)$ and 
	an unique coninuously differentiable function 
	$v:U \to \R^{n_y}$ with $v(t_0,x_0,u_0) = v_0$ and $\tilde{g}(t,x,u,v(t,x,u)) = 0$ for all $(t,x,u) \in U$.\\
	Let $x,\, u$ be the locally unique solution of the coupled initial value problems 
	\begin{subequations}\label{eq:help_prob}
	\begin{align}
		x^{\prime}(t) &= f(t,x(t),u(t)+Qv(t,x(t),u(t))),\label{eq:help_prob_1}\\
		u^{\prime}(t) &= Pv(t,x(t),u(t)),\label{eq:help_prob_2}\\
		x(t_0) &= x_0 \label{eq:help_prob_3}\\
		u(t_0) &= Py_0\label{eq:help_prob_4}
	\end{align}
	\end{subequations}
	with $ t \in \tilde{I}\subset I$ an open interval with $t_0 \in \tilde{I}$.
	
	We set $y(t) := u(t) + Qv(t,x(t),u(t))$, it remains to show, that $(x,y)$ is a solution to the reduced DAE \eqref{eq:reducesem_b}. 
	Therefore we first prove, that it holds $Qu(t) = 0$ for all $t \in \tilde{I}$. 
	Multiply \eqref{eq:help_prob_2} and \eqref{eq:help_prob_4} by Q implies 
	\begin{subequations}
	\begin{align}
		(Qu)^{\prime}(t) &= Qu^{\prime}(t) = QPv(t,x(t),u(t)) = 0,\\
		Qu(t_0) &= QPy_0 = 0.
	\end{align}
	\end{subequations}
	Then for $\alpha(t) = Qu(t)$ it follows $\alpha^{\prime}(t) = 0$ and $\alpha(t_0) = 0$, which implies $\alpha(t) = 0$ for all $t \in \tilde{I}$ 
	and $Qu(t) = 0$. Especially we have $Pu(t) = u(t)$ for all $t \in \tilde{I}$. 
	For $x$ and $y$ with \eqref{eq:eq_second}, \eqref{eq:eq_third} and \eqref{eq:help_prob_2} it follows
	\begin{align*}
		x^{\prime}(t) &= f(t,x(t),y(t)),\\
		0 &= \tilde{g}(t,x(t),u(t),v(t,x(t),u(t)))\\
			&=h(t,x(t),y(t))+\partial_yg(t,x(t),y(t))v(t,x(t),u(t))\\
			&=h(t,x(t),y(t))+\partial_yg(t,x(t),y(t))Pv(t,x(t),u(t))\\
			&=h(t,x(t),y(t))+\partial_yg(t,x(t),y(t))u^{\prime}(t)\\
			&=h(t,x(t),y(t))+\partial_yg(t,x(t),y(t))(Pu^{\prime}(t)+0)\\
			&=h(t,x(t),y(t))+\partial_yg(t,x(t),y(t))(Pu^{\prime}(t)+PQ\partial_tv(t,x(t),u(t))\\
			&=h(t,x(t),y(t))+\partial_yg(t,x(t),y(t))Py^{\prime}(t)\\
			&=h(t,x(t),y(t))+\partial_yg(t,x(t),y(t))y^{\prime}(t)\\
			&=\partial_tg(t,x(t),y(t))+\partial_xg(t,x(t),y(t))f(t,x(t),y(t))+\partial_yg(t,x(t),y(t))y^{\prime}(t).
	\end{align*}
	and
	\begin{align*}
	 x(t_0) = x_0, \, y(t_0) = u(t_0) + Qv(t_0,x(t_0),u(t_0)) = u_0 + Qv_0 = Py_0 + Q(Pw_0+Qy_0) = Py_0 + Qy_0 = y_0.
	\end{align*}		 
	Therefore we have a locally unique solution of the reduced problem for the given initial values and the claim follows.\hfill$\Box$
\end{proof}

\section{Mathematical Representation of District Heating Networks}
\label{sec:3}

In the following we give a short description of the modeling of a district heating network. For more details we refer the reader to \cite{Borsche19, Koecher00, Krug19}, for instance. The goal is to put the representation of a district heating network into a differential algebraic system of equations in a semi-explicit form like \eqref{eq:semiDAE}. For the modeling of the network we refer to \cite{Herty2006, Qiu18} for a graph theoretical background we refer the reader to \cite{BangGutinGregory2008, Diestel2017}. 

\subsection{Modeling of district heating networks}
\subsubsection{A single pipe}
\label{sec:2.1}
A network of a district heating system is composed by several different components. First we describe the mathematical
model of the flow in a single pipe. These can be joined to build a network and connected to households via suitable
coupling conditions. We follow \cite{Borsche19} and consider the following three equations
\begin{subequations}\label{eq:afangg}
\begin{align}
\label{eq:anfanggl1}
	\partial_t\rho + \partial_x q &= 0\\
\label{eq:anfanggl2}
	\partial_t q + \partial_x \left( \frac{q^2}{\rho}+p \right) &= -\frac{\lambda}{2d}\frac{q\vert q \vert}{\rho} - g(\partial_x h)\rho \\
\label{eq:anfanggl3}
	\partial_t(c_p\rho T) + \partial_x(qc_pT) &= - \frac{4k}{d}(T-T_{ext})
\end{align}
\end{subequations}
which describe the conservation of mass, the balance of momentum and energy, respectively. 
Here $\rho$ denotes the density of the water, $q = \rho v$, where $v$ is velocity, $p$ is the pressure and $T$ is the
temperature of the fluid. The parameter $\lambda$ is a friction coefficient for the Darcy-Weisbach friction formula and $d$ the diameter of the pipe.
The parameter $c_p$ is the specific heat
capacity of the water. The term $g(\partial_x h)$ takes the vertical displacement $h$ into account, where $g$ is the gravitational acceleration. In a pipe with length $L$ and height difference
$\Delta h $ it holds $\partial_x h = \frac{\Delta h}{L}$. Therefore we set $g(\partial_x h) = g \Delta h/L$ 
where $g\approx 9.80665\, ms^{-2}$ is the gravitational acceleration and $L$ the length of the pipe.
The right-hand side in equation \eqref{eq:anfanggl3} models the cooling related to the outer temperature
$T_{ext}$ with the thermal transmittance $k$.
As the water in the pipes is almost incompressible and the temperature difference is not that big we assume that the density $\rho$ is constant so that we have $\partial_t\rho = 0$. 
Inserting this assumption in equation \eqref{eq:anfanggl1} it follows 
\begin{equation}
	\label{vkonstant}
	\partial_x v = 0.
\end{equation}
Using $\partial_t\rho = 0$ and inserting \eqref{vkonstant} in \eqref{eq:anfanggl3} we obtain
\begin{align*}
	c_p\rho\partial_tT + c_p\rho v\partial_xT &= - \frac{4k}{d}(T-T_{ext}).
\end{align*} 
The fact that $v$ is constant in space transforms equation \eqref{eq:anfanggl2} to the incompressible Euler equation
\begin{align}\label{incompeuler}
	\rho \partial_t v + \partial_x p = -\frac{\lambda}{2d}v\vert v \vert \rho - g(\partial_x h)\rho.
\end{align}
In the next step we integrate \eqref{incompeuler} over the length of the pipe $[0,L]$ and get
\begin{align*}
	\partial_t v\int_0^L \rho \,\mathrm dx +\int_0^L\partial_x p \,\mathrm dx &=-\frac{\lambda}{2d}v\vert v \vert\int_0^L \rho \,\mathrm dx - g(\partial_x h) \int_0^L \rho\,\mathrm dx,\\
	\rho\partial_t v +\frac{p(t,L)-p(t,0)}{L} &=-\frac{\lambda}{2d}v\vert v \vert \rho - g(\partial_x h) \rho.
\end{align*}
Moreover, we assume that the sign of the velocity is always positive. 
Therefore, we can simplify $v\vert v \vert$ to $v^2$. 
In summary, we get for a single pipe the two equations
\begin{subequations}\label{eq:endeq}
\begin{align}
\label{eq:transport}
\partial_t T(t,x) + v(t) \partial_x T(t,x) &= - \frac{4k}{c_p d \rho} (T(t,x) - T_{ext}),\\
\label{eq:euler}
\rho \partial_t v(t) + \frac{p(t,L) - p(t,0)}{L} &= - \frac{\lambda}{2d} v(t)^2 \rho - g(\partial_x h) \rho.
\end{align}
\end{subequations}
At the moment the transport equation is a partial differential equation, which leads to a system of partial differential algebraic equations (PDAEs), we are 
only working with DAEs at the moment, therefore we apply a spatial discretization (implicit Euler) by the methods of lines for the transport equation \eqref{eq:transport} and get 
\begin{align}\label{eq:transportdisc}
	\partial_t T_j(t) + \frac{v(t)}{\Delta x}(T_j(t)-T_{j-1}(t)) &= - \frac{4k}{c_p d \rho} (T_j(t) - T_{ext}) \mbox{ for } j = 2,...,n_1.
\end{align}
For a list of all parameters and variables of the model see Table \ref{tab:Var}.  
The friction factor $\lambda$ for a turbulent flow is modeled by the flow-independent law of Nikuradse (see e.g., \cite{Fuegenschuh15}), i.e.,
\begin{align*}
\lambda = \left( 2 \log_{10} \left(\frac{d}{k_{rough}} \right) + 1.138 \right)^{-2},
\end{align*}
where $k_{rough}$ is the roughness of the inner pipe wall.
Notice, in a single pipe the system composed by \eqref{eq:endeq} has to be supplemented with boundary conditions and initial datum.
  A common choice for $v>0$ is to provide $p$ and $T$ at the left (the power plant), and specify the demanded flow $q$ at 
 the right end (consumer site). For a pipe in a network these boundary conditions are replaced by coupling conditions in the junctions. 

\begin{table}[!t]
\caption{Variables (top) and parameters (bottom) of the district heating network model.}
\label{tab:Var}       
\begin{tabular}{llcc}
\toprule
Symbol & Explanation & Unit & Example    \\
\midrule
$v(t)$ & Flow velocity  & m s$^{-1}$ & -- \\
$p(t,0)$, $p(t,L)$ & Pressure at the left, right end of a pipe & Pa$=$kg m$^{-1}$s$^{-2}$& \phantom{1}$5\cdot 10^5, 2\cdot 10^5$\\
$T(t,x)$  & Water temperature in a pipe  & C& -- \\
$q(t)$ & Mass flow in a pipe, $q(t)=\rho v(t)$  & kg m$^{-2}$s$^{-1}$&  -- \\
$T_{in}(t)$ & Temperature of hot water entering the network & C & -- \\
$p_{in}(t)$ & Pressure of the entering hot water  & \phantom{12}Pa$=$kg m$^{-1}$s$^{-2}$ & -- \\
\midrule
$\rho$ & Density of the water & kg m$^{-3}$& $960$ \\
$t$ & Time coordinate; $t \in \mathcal{T}$  & s& --\\
$\mathcal{T}$ & Time horizon; $\mathcal{T} = [t_0,t_f]$  & -- &$[0,4000]$\\
$x$ & Spatial coordinate in a pipe & m &--\\
$L$ & Length of a pipe & m &$300$\\
$\Delta h$ & Height difference in a pipe & m & -- \\
$d$ & Diameter of a pipe & m &$0.1$\\
$A$ & Cross-sectional area of a pipe; $A=\pi(\frac{d}{2})^{2}$ & m$^{2}$ & $7.85\cdot 10^{-3}$\\
$\lambda$ & Friction factor of a pipe & $1$&$0.00251$ \\
$Q_k(t)$ & Power consumption of a consumer & W & $1\cdot 10^6$\\
$k_{rough}$ & Roghness of the inner wall of a pipe a & m &$0.00026$\\
$k$ & Heat transfer coefficient of the wall of a pipe & W m $^{-2}$ C $^{-1}$ &$0.1$\\
$T_{out}$ & Consumers' outlet water temperature & C& $60$\\
$T_{ext}$ & Surrounding temperature & C& $20$\\
$c_p$ & Specific heat capacity of water & J kg$^{-1}$C$^{-1}$ &$4160$\\
$g(\partial_x b)$ & Gravitational acceleration $=g \cdot \frac{\Delta h}{L} $  & m s$^{-2}$ &$9.80665$\\
\bottomrule
\end{tabular}
\end{table}

\subsubsection{Network}
\label{subsec:2.1.2}

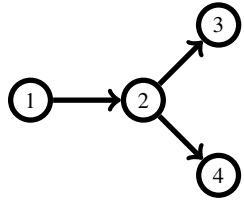
\begin{figure}[h]
\begin{center}
\begin{tikzpicture}[scale=1]
	\node[line width=2pt] (v1) at (0,1) [circle,draw] {$1$};
	\node[line width=2pt] (v2) at (1.5,1) [circle,draw] {$2$};
	\node[line width=2pt] (v3) at (2.5,2) [circle,draw] {$3$};
	\node[line width=2pt] (v4) at (2.5,0) [circle,draw] {$4$};
	
	\draw[->][line width=2pt] (v1) to node[above] {} (v2);
	\draw[->][line width=2pt] (v2) to node[above] {} (v3);
	\draw[->][line width=2pt] (v2) to node[above] {} (v4);
\end{tikzpicture}
\end{center}
\caption{Illustration of a junction.}
\label{fig:junc}   
\end{figure}

In a district heating network the hot water is distributed to households via a system of pipes. A network of the same structure 
transports the colder water back to the power plant. \\
Such networks can be modeled by prescribing suitable coupling conditions at the junctions additional to above equations on the edges
 (compare to Fig \ref{fig:junc}). 
 In both networks we consider the following coupling conditions in every interior node $j$:
\begin{subequations}
\label{couplingCond}
\begin{align}
\label{conserveMass}
\sum_{e_i \in \sigma_j} A_i q_i &= \sum_{i\in \Sigma_j}A_iq_i,  \\
\label{conserveEnergy}
\sum_{e_i \in \sigma_j} c_p A_i q_i T_i(t,L_i) &= \sum_{i\in \Sigma_j}c_p A_i q_i T_i(t,0), \\
\label{contPressure}
p_i(t,L_i) &= p_l(t,0) && \text{for all } i \in \sigma_j, \, l\in \Sigma_j, \\
\label{perfectMix}
T_i(t,0) &= T_l(t,0) && \text{for } i,l \in \Sigma,\, i \neq l,
\end{align}
\end{subequations}
where $\sigma_j$ is the set of all pipes incoming pipe $j$ and $\Sigma_j$ is the set of all pipes leaving pipe $j$. 
The junction (Fig \ref{fig:junc} at $2$) is assumed to connect $|\sigma_j|+|\Sigma_j|$ pipes and $A_i$ 
denotes the cross section of the $i$th pipe. Equation \eqref{conserveMass} states the conservation of mass and 
\eqref{conserveEnergy} the conservation of energy. The continuity of the pressure \eqref{contPressure} is a widely 
used condition, see e.g. \cite{Banda06, Colombo08, Domschke15}. Additionally we assume a perfect mixing of flows at the 
junction, which means that we assume the same temperature in all outgoing pipes, \eqref{perfectMix}. \\

An additional important component in a district heating network are the consumers. Each consumer is demanding a certain
 amount of thermal power $Q_k(t)$. Further the outgoing temperature $T_{out}$ is assumed to be a fixed value and no mass is lost. 
 This leads to the following equations:
\begin{subequations}
\label{connectNetwork}
\begin{align}
\label{noMassLost}
q_{in} &= q_{out}, \\
\label{consumerDemand}
Q_k(t) &= c_p A q_{in}(T_{in} - T_{out}),
\end{align}
\end{subequations}
where $T_{in}$ is the temperature of the flow arriving at the household. These relations of the three quantities $q_{in}$, $q_{out}$ 
and $T_{in}$ are connecting the supplying network with the one for the return flow.
Finally we need to introduce the power plant. We will assume that for the simulation a given temperature field will enter the network from the power plant and a pressure field will leave the 
network to the power plant. 
\begin{subequations}\label{eq:simulate}
\begin{align}
\label{eq:simulate1}
p_N(t,L_N) &= p_{in}(t), \\
T_1(t,0) &= T_{in}(t)\label{eq:simulate2}
\end{align}
\end{subequations}
for some functions $p_{in}$ and $T_{in}$.
\subsection{Graph theoretical modeling}
\label{sec:2.2}
The abstract network is described by a directed graph 
\begin{equation}\label{graph}
	\g = ( \n,\e).
\end{equation}
Here $\n$ denotes the set of nodes, which consist of the set of supply nodes $\n_s$, demand nodes $\n_d$ and interior nodes $\n_0$ of 
the network. Here, the supply nodes represent the set of nodes in the network, where water
is injected into the network. The demand nodes form a set of nodes, where the water is extracted from the network and interior nodes are the rest. 
Sometimes interior nodes are called junction nodes. 
We assume from now on that demand nodes and supply nodes are the 
only boundary nodes. That means they are only connected to one pipe. If supply or demand nodes exist that are connected to more than one pipe, 
we add a short pipe to that node and declare the new node as the demand 
or the supply node and the old one becomes an interior node.
This short pipe is sometimes called pseudo or virtual pipe. We consider only networks with a tree configuration and with a single power plant, compare e.g. \cite{bordin2016optimization}. In this case, 
the flow direction is a-priori defined as the water flows from the power plant to the consumers and there are no loops in the system.
\begin{definition}
	The nodes of a graph $\g$, which connected at least two pipes, are called junction nodes. 
\end{definition}
\begin{definition}
	A polytree is a directed acyclic graph whose underlying undirected graph is a tree.
\end{definition}
Let $n_s $ be the number of supply nodes, $n_d$ is the number of demand nodes, $n_{junc}$ is the number of interior nodes and 
$\oN := n_s+n_d+n_{junc}$ is the number of all nodes. The set $\e \subset \n \times \n$ contains the pipes of the network in the sense, that $e = (v_1,v_2) \in \e$ describes the pipe 
between node $v_1$ and $v_2$ with the water direction from $v_1$ to $v_2$.  
A pipe attached to a supply node is called a supply pipe, while a pipe attached to a demand pipe is called a demand pipe. 
A supply pipe is directed away from the supply node and a demand pipe is directed towards the demand 
node. We assume that $N$ is the number of pipes.
To model the network 
with consumers and a power plant, we assume that a consumer or a power plant is a ''break'' in the network in the sense, 
that a consumer or a power plant is localized between a demand and supply node. 
For every supply node the temperature during the process is known.\\
To work with the graph, we number the pipes and the nodes. This numberation implies an order in the network. Of course, there is more than one possiblity, to specify the order we use the 
following lemma.
\begin{lemma}
	Given a directed acyclic graph with no loops, we can order the pipes in such a way that at every node all incoming pipes have a 
	lower order as all the outgoing pipes, or it does not have an incoming pipe. We call this ordering \textit{direction
	following ordering}. 
\end{lemma}
\begin{proof}
We skip a detailed description of the proof here and refer the reader to \cite[Section 4.2]{Qiu18}.\hfill$\Box$
\end{proof}
We assume that we have $n_c \in \mathbb{N}$ consumers in the network. From now, the first $n_s-n_c$ supply nodes are incoming nodes after a power plant. 
Also the last $n_d-n_c$ demand nodes are outgoing nodes before a power plant. We also assume that the last $n_d-n_c$ pipes are demand pipes before a power plant. Note 
that the direction following ordering is not unique. Since we assume that we only have one power plant in the system it holds $n_s -n_c = n_d - n_c = 1$. \\
The incidence matrix shows the relationship between the nodes $\n$ and pipes $\e$. 
The matrix has one column for each pipe $e_i \in \e$ and one row for each node $\on_j \in \n$.  
\begin{definition}
	The \emph{incidence matrix} of a directed graph $\g = (\n,\e)$ is a $\oN \times N$ matrix $\mathcal{I}$ such that 
	\begin{equation*}
		\mathcal{I}_{j,i} := \begin{cases} 1 & \mbox{ if pipe } e_i \mbox{ leaves node } \on_j \\
					-1 & \mbox{ if pipe } e_i \mbox{ enters node } \on_j \\
					0 & \mbox{otherwise} \end{cases}.
	\end{equation*}
\end{definition}
To describe the consumers in the network we define so called \emph{consumer matrices} $\c^1$ and $\c^2$ with 
\begin{equation*}
	\c^1_{k,i} =  \begin{cases} 1 & \mbox{ if consumer } k \mbox{ is localized after pipe } e_i \\
						0 & \mbox{otherwise} 
			\end{cases} 
\end{equation*}
and
\begin{equation*}
	\c^2_{k,i} =  \begin{cases} 1 & \mbox{ if consumer } k \mbox{ is localized before pipe } e_i \\
						0 & \mbox{otherwise} 
			\end{cases}.
\end{equation*} 
Here and in what follows we set
\begin{align*}
	x(t) &= \begin{pmatrix}
		x_1(t) \\
		x_2 (t)
	\end{pmatrix},\quad
	x_2(t) = \begin{pmatrix}
		v_1(t)\\
		\vdots\\
		v_N(t)
	\end{pmatrix} \in \R^N,  \\
	x_1(t) &= \begin{pmatrix} 
		T_{1,2}(t) \\ 
		\vdots \\ 
		T_{1,n_1}(t) \\ 
		\vdots \\ 
		T_{N,2}(t) \\ 
		\vdots \\ 
		T_{N,n_N}(t) 
	\end{pmatrix} \in \R^{\tilde{n}}, \quad 
	y(t) = \begin{pmatrix}
		T_{1,1}(t) \\ 
		\vdots \\ 
		T_{N,1}(t) \\ 
		p_1(t,0) \\ 
		p_1(t,L_1) \\ 
		\vdots \\ 
		p_N(t,L_N) 
	\end{pmatrix} \in \R^{3N}
\end{align*}
with $n = \sum_{i=1}^N n_i, \, \tilde{n} = n-N$.
Now a semi-discrete version of the network model discribed by \eqref{eq:euler}, \eqref{eq:transportdisc}, \eqref{couplingCond}, \eqref{connectNetwork} and \eqref{eq:simulate} is 
given by the following semi-explicit DAE
\begin{align}
\label{eq:semi_expl}
\begin{split}
\dot{x}_1(t) &= f_1(t, x_1(t), x_2(t), y(t)) \\
\dot{x}_2(t) &= f_2(t, x_1(t), x_2(t), y(t)) \\
0 &= g_1(t, x_1(t), x_2(t), y(t)) \\
0 &= g_2(t, x_1(t), x_2(t))
\end{split}
\end{align}
where $f_1: \R \times \R^{n_{x_1}} \times \R^{n_{x_2}} \times \R^{n_y} \rightarrow \times \R^{n_{x_1}}$, 
$f_2: \R \times \R^{n_{x_1}} \times \R^{n_{x_2}} \times \R^{n_y} \rightarrow \times \R^{n_{x_2}}$, 
$g_1 : \R \times \R^{n_x} \times \R^{n_y} \rightarrow \R^{n_{y_1}}$ and 
$g_2 : \R \times \R^{n_x} \rightarrow \R^{n_{y_2}}$. Note, that $n_{x_1} + n_{x_2} = n_x$ and that $n_{x_1} = \tilde{n}$, $n_{x_2} = N$.
For a detailed description of $f_1,$ $f_2$, $g_1$ and $g_2$ see Appendix A. 
\begin{lemma}
\label{lem:index}
Assume that district heating is modeled by a polytree $\g$ and every demand node is a outgoing node localized in front of a consumer or a power plant. For all consumers it holds that
\begin{align}\label{assumptionQpos1}
Q_k(t) > 0 \mbox{ for all }t\in I, \qquad \text{k=1,\ldots, $n_c$}.
\end{align}
Then the semi-explicit DAE \eqref{eq:semi_expl} has differentiation index $2$ and it holds $v(t) >0$ for all $t \in I$ in every pipe. 
\end{lemma}
\begin{proof}
The consumer demand equation \eqref{consumerDemand}, assumption \eqref{assumptionQpos1} and the fact that 
every demand node is localized in front of a consumer imply for physical reasons that the velocity $v(t)$ 
in every demand pipe is strictly positive for all $t \in I$. Due to the tree structure, the velocity $v(t)$
 in every pipe need to be strictly positive for all $t \in I$. \\
Similar calculations like in \eqref{calcIndex} show that the equations in $g_1$ have differentiation index $1$. 
For $g_2$ we use the same argument two times 
to  show that the equations have differentiation index $2$. \hfill$\Box$ 
\end{proof}	
To reduce the index, one can use the method discribed in Section \ref{sec:intro}.
Another way occurs,  if we neglect the term $\dot{v}$ in the modified incompressible Euler equation \eqref{eq:euler}. 
This implies the new semi-explicit DAE
\begin{align}
\label{eq:semi_expl_index1}
\begin{split}
x_1'(t) &= f_1(t, x_1(t), x_2(t), y(t)) \\
0 &= f_2(t, x_1(t), x_2(t), y(t)) \\
0 &= g_1(t, x_1(t), x_2(t), y(t)) \\
0 &= g_2(t, x_1(t), x_2(t)),
\end{split}
\end{align}
Since we work with water in the district heating network, the pressure differences in the pipes and the velocity changes are restricted through the natural behavior of water. 
Therefore the assumption $\vert \dot{v} \vert \approx 0$ does not produce a significant error in the solution of the original DAE \eqref{eq:semi_expl} and the new semi-explicit DAE \eqref{eq:semi_expl_index1}. 
\begin{lemma}
\label{lem:index1}
Assume that a district heating is modeled with a polytree $\g$ and every demand node is a outgoing node localized in front of a consumer or a power plant and for all consumers it holds that
\begin{align*}
Q_k(t) > 0\quad\text{for all }t\in I\text{ and }k=1,\ldots, n_c.
\end{align*}
Then the semi-explicit DAE \eqref{eq:semi_expl_index1} has differentiation index $1$ and it holds $v(t) >0$ for all $t \in I$ in every pipe. 
\end{lemma}
\begin{proof}
The proof utilizes the same arguments as in the proof of Lemma \ref{lem:index}.\hfill$\Box$
\end{proof}
\begin{definition}
We call the semi-explicit DAE \eqref{eq:semi_expl} full problem or index $2$ problem and the semi-explicit DAE \eqref{eq:semi_expl_index1} reduced prolbem or index $1$ problem. 
\end{definition}

\subsection{Examples}
\label{subsec:Examples}	
\begin{figure}[h]
\begin{subfigure}{.5\linewidth}
\centering
\begin{tikzpicture}[scale=1.5]
\node[line width=2pt,black] (v1) at (0,1) [circle,draw] {$1$};
\node[line width=2pt,black] (v2) at (2,0) [circle,draw] {$2$};
\node[line width=2pt,black] (v3) at (2.5,-1) [circle,draw] {$3$};
\node[line width=2pt,black] (v4) at (1,1) [circle,draw] {$4$};
\node[line width=2pt,black] (v5) at (1,-0.5) [circle,draw] {$5$};
\node[line width=2pt,black] (v6) at (2,0.5) [circle,draw] {$6$};
\node[line width=2pt,black] (v7) at (2.5,1.5) [circle,draw] {$7$};
\node[line width=2pt,black] (v8) at (0,-0.5) [circle,draw] {$8$};
\draw[->, line width=2pt,] (v1) to node[above] {$1$} (v4);
\draw[->, line width=2pt,] (v2) to node[above] {$2$} (v5);
\draw[->, line width=2pt,] (v3) to node[above] {$3$} (v5);
\draw[->, line width=2pt,] (v4) to node[above] {$4$} (v6);
\draw[->, line width=2pt,] (v4) to node[above] {$5$} (v7);
\draw[->, line width=2pt,] (v5) to node[above] {$6$} (v8);
\draw[-, line width=2pt,gray] (v6) to node[left] {Cons. 1} (v2);
\draw[-, line width=2pt,gray] (v7) to node[right] {Cons. 2} (v3);
\draw[-, line width=2pt,gray] (v1) to node[right] {Power 1} (v8);
\end{tikzpicture}
\caption{Network modeling}
\label{fig:example1a}      
\end{subfigure}
\begin{subfigure}{.5\linewidth}
\centering
\begin{tikzpicture}[scale=0.75]
 Knoten
\node[line width=2pt,black] (v1) at (0,0) [circle,draw] {$1$};
\node[line width=2pt,black] (v2) at (1,0) [circle,draw] {$2$};
\node[line width=2pt,black] (v3) at (2,0) [circle,draw] {$3$};
\node[line width=2pt,black] (v4) at (3,0) [circle,draw] {$4$};
\node[line width=2pt,black] (v5) at (4,0) [circle,draw] {$5$};
\node[line width=2pt,black] (v6) at (5,0) [circle,draw] {$6$};
\node[line width=2pt,black] (v7) at (6,0) [circle,draw] {$7$};
\node[line width=2pt,black] (v8) at (7,0) [circle,draw] {$8$};
%
\draw[->, line width=2pt] (v1) to [bend right=60] node[below] {$1$} (v4);
\draw[->, line width=2pt] (v4) to [bend right=40] node[below] {$4$} (v6);
\draw[->, line width=2pt] (v4) to [bend right=60] node[below] {$5$} (v7);
\draw[->, line width=2pt] (v2) to [bend left =60] node[above] {$2$} (v5);
\draw[->, line width=2pt] (v3) to [bend left =40] node[above] {$3$} (v5);
\draw[->, line width=2pt] (v5) to [bend left =60] node[above] {$6$} (v8);
\end{tikzpicture}
\caption{Direction following ordering}
\label{fig:example1b}      
\end{subfigure}
\caption{Example of a network with two consumers}
\label{fig:example33}
\end{figure}
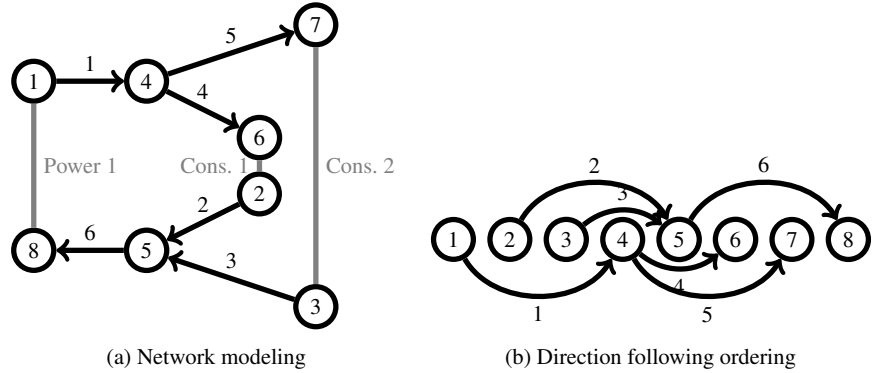

For the numerical tests we introduce one example with two consumers. In figure \ref{fig:example1a} the network model for the 
example with two consumers is shown and figure \ref{fig:example1b} shows the direction following ordering for the network.

\noindent The incidence matrix for the example has the following structure:  
\begin{align*}
\mathcal{I} = \begin{pmatrix}
	1 & 0 & 0 & 0 & 0 & 0 \\
	0 & 1 & 0 & 0 & 0 & 0 \\
	0 & 0 & 1 & 0 & 0 & 0 \\
	-1& 0 & 0 & 1 & 1 & 0 \\
	0 &-1 &-1 & 0 & 0 & 1 \\
	0 & 0 & 0 & -1& 0 & 0 \\
	0 & 0 & 0 & 0 &-1 & 0 \\
	0 & 0 & 0 & 0 & 0 &-1 \\
\end{pmatrix}.
\end{align*}
The consumer matrices are
\begin{align*}
\c^1 = \begin{pmatrix}
	0 & 0 & 0 & 1 & 0 & 0 \\
	0 & 0 & 0 & 0 & 1 & 0 
\end{pmatrix} \mbox{ and }
\c^2 = \begin{pmatrix}
	0 & 1 & 0 & 0 & 0 & 0 \\
	0 & 0 & 1 & 0 & 0 & 0 
\end{pmatrix}.
\end{align*}
It holds $N = 6,\, \oN = 8,\, n_s = 3,\, n_d = 3,\, n_{junc} = 2,\, n_c=2$.

\section{Numerical Experiments}
\label{sec:4}

In this section we present our numerical experiments. Here and in the following we will use IDA from SUNDIALS (see: \cite{Sundials}) as DAE solver. IDA is implemented for python in e.g. the scikit package odes, see \cite{odes} and assimulo, see \cite{Andersson2015}.

In the following examples it is necessary to calculate consistent initial values (see Definition \ref{def:consistent_initial_value}). For this we solve a constrained least-squares problem compare e.g. \cite{burger2017survey}. More precisely, let $z_0=(\dot{x}_0, x_0, y_0)$ be the initial value we are looking for. Then we solve
\begin{align*}
\min \frac{1}{2} \Vert z - z_0 \Vert^2_2
\end{align*}
for a fixed $z$ which is close to the searched initial value, subject to the constraints
\begin{align*}
0 = \begin{pmatrix}
x_0' - f(t_0, x_0, y_0) \\
g_1(t_0, x_0, y_0) \\
g_2(t_0, x_0, y_0)
\end{pmatrix}
\end{align*}
and in the index $2$ case we need in addition the hidden constraint
\begin{align*}
0 = \frac{\mathrm d}{\mathrm dt}g_2(t_0, x_0, y_0).
\end{align*}
This problem is then finally solved with a trust-region SQP procedure.

\subsection{Compare different spatial discretization}
\label{sec:4.1}
Before we start some numerical experiments we want to validate our implementation for a special case, where we know the analytic solution. We will have a look on the behavior for different spatial discretizations as discussed later in Remark \ref{rem:spatial_discretization}. Also we will compare the behavior of the full and the reduced problem, compare Lemma \ref{lem:index}. In both cases we compare different spatial discretizations, namely a first-, second- and third-order discretization, see Appendix A. We will use the network from the example in Section \ref{subsec:Examples}, see Figure \ref{fig:example33}. In both cases we use the parameters
\begin{align*}
k &= -1,    & c_p     &= 2, & d &= 1, & \rho            &= 2, \\
T_{ext}&=0, & \lambda &= 2, & L &= 1, & g(\partial_x h) &= 1,
\end{align*}
the consumer demands
\begin{align*}
Q_1(t) &= \frac{1}{3} (2 \exp(\frac{3}{2})-1) \pi \exp(1+t), \\ 
Q_2(t) &= \frac{1}{6} (2 \exp(3)-1) \pi \exp(1+t)
\end{align*}
and for the velocities and the temperatures the functions 
\begin{align*}
v_1(t)&= \frac{1}{2-t}, & T_1(t,x) &= \exp(t+x)(2-t), \\
v_2(t)&= \frac{2}{6-3t}, & T_2(t,x) &= \frac{1}{2}\exp(1+t+\frac{3}{2}x)(2-t), \\
v_3(t)&= \frac{1}{6-3t}, & T_3(t,x) &= \frac{1}{2}\exp(1+t+3x)(2-t), \\
v_4(t)&= \frac{2}{6-3t}, & T_4(t,x) &= \exp(1+t+\frac{3}{2}x)(2-t), \\
v_5(t)&= \frac{1}{6-3t}, & T_5(t,x) &= \exp(1+t+3x)(2-t), \\
v_6(t)&= \frac{1}{2-t}, & T_6(t,x) &= \frac{1}{6}(2+\exp(\frac{3}{2}))\exp(\frac{5}{2}+t+x)(2-t) + 5.
\end{align*}
In the index $1$ case, where we neglect the term $\dot{v}$, we set for the pressures
\begin{align*}
p_1(t,0) &= \frac{3}{(t-2)^2} + 2, & p_4(t,0) &= \frac{1}{(t-2)^2}, \\
p_1(t,L) &= \frac{1}{(t-2)^2}, & p_4(t,L) &= \frac{1}{9(t-2)^2}-2, \\
p_2(t,0) &= \frac{44}{9(t-2)^2} + 4, & p_5(t,0) &= \frac{1}{(t-2)^2}, \\
p_2(t,L) &= \frac{4}{(t-2)^2} + 2, & p_5(t,L) &= \frac{7}{9(t-2)^2}-2, \\
p_3(t,0) &= \frac{38}{9(t-2)^2} + 4, & p_6(t,0) &= \frac{4}{(t-2)^2}+2, \\
p_3(t,L) &= \frac{4}{(t-2)^2} + 2, & p_6(t,L) &= \frac{2}{(t-2)^2} 
\end{align*}
and in the index $2$ case, where we do not neglect the term $\dot{v}$, we set
\begin{align*}
p_1(t,0) &= \frac{5}{(t-2)^2} + 2, & p_4(t,0) &= \frac{1}{(t-2)^2}, \\
p_1(t,L) &= \frac{1}{(t-2)^2}, & p_4(t,L) &= -\frac{11}{9(t-2)^2}-2, \\
p_2(t,0) &= \frac{74}{9(t-2)^2} + 4, & p_5(t,0) &= \frac{1}{(t-2)^2}, \\
p_2(t,L) &= \frac{6}{(t-2)^2} + 2, & p_5(t,L) &= \frac{1}{9(t-2)^2}-2, \\
p_3(t,0) &= \frac{62}{9(t-2)^2} + 4, & p_6(t,0) &= \frac{6}{(t-2)^2}+2, \\
p_3(t,L) &= \frac{6}{(t-2)^2} + 2, & p_6(t,L) &= \frac{2}{(t-2)^2}.
\end{align*}
\begin{figure}
\begin{subfigure}{.5\linewidth}
	\centering
	\includegraphics[width=\textwidth]{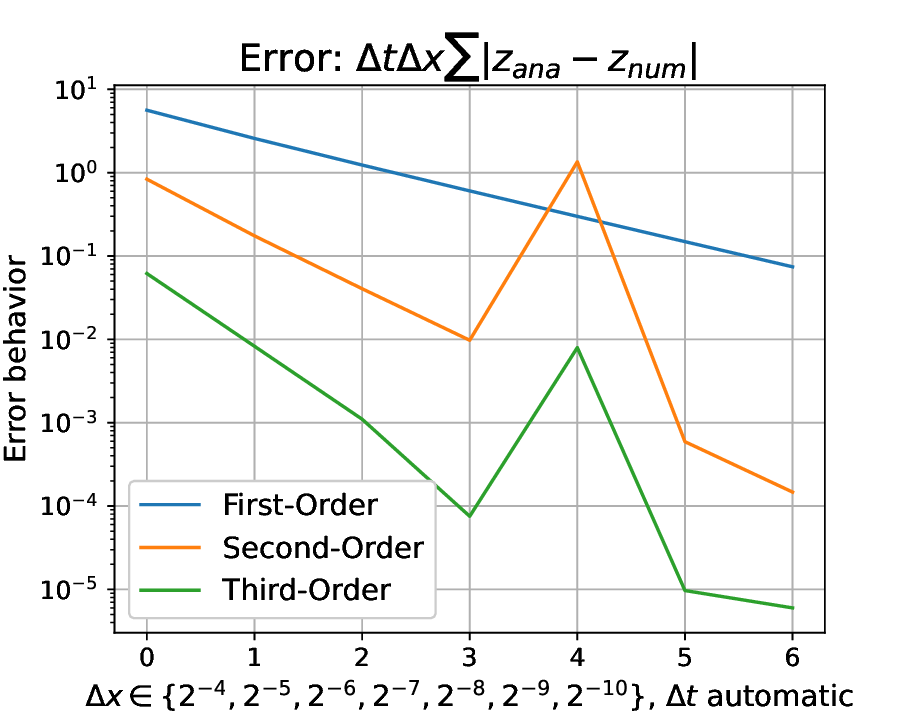}
	\caption{discete $L_1$-error}
	\label{fig:TwoConsumer_index1_abs}
\end{subfigure}
\begin{subfigure}{.5\linewidth}
	\centering
	\includegraphics[width=\textwidth]{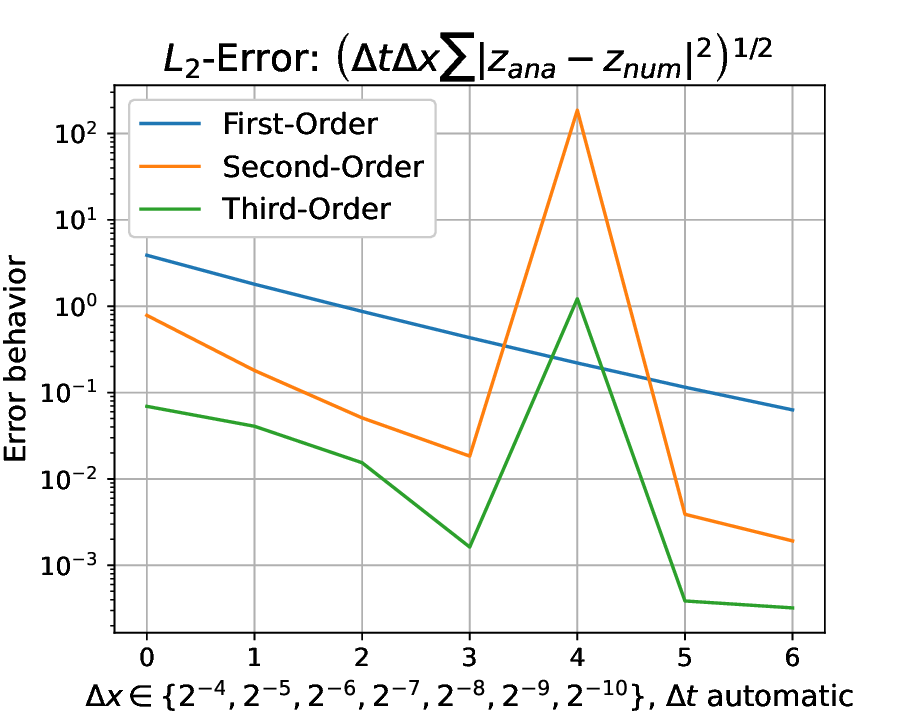}
	\caption{discrete $L_2$-error}
	\label{fig:TwoConsumer_index1_l2}
\end{subfigure}
\begin{subfigure}{.5\linewidth}
	\centering
	\includegraphics[width=\textwidth]{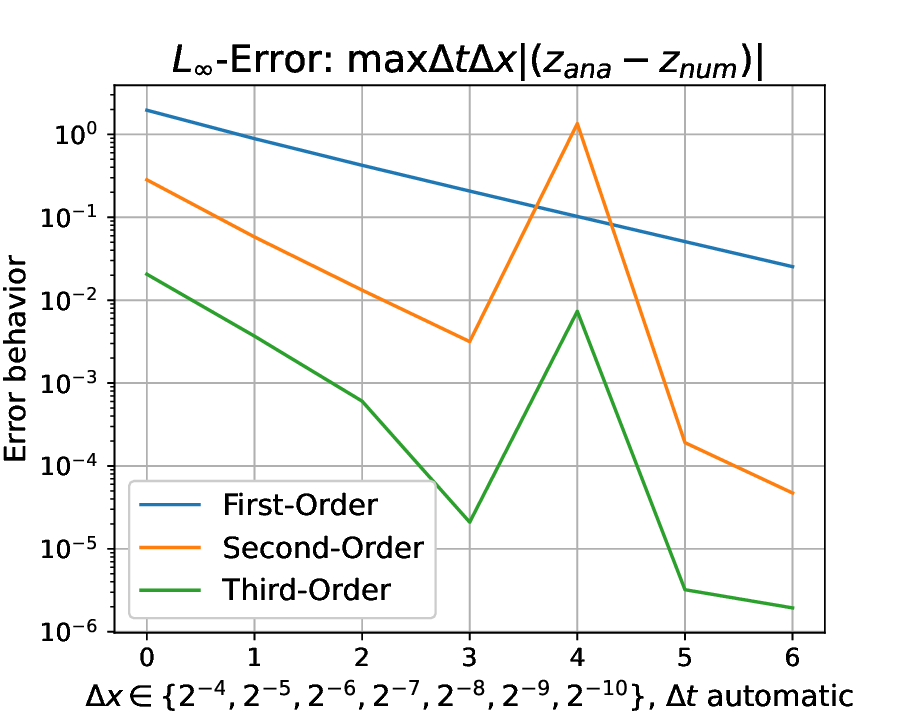}
	\caption{discrete $L_{\infty}$-error}
	\label{fig:TwoConsumer_index1_linf}
\end{subfigure}
\begin{subfigure}{.5\linewidth}
	\centering
	\includegraphics[width=\textwidth]{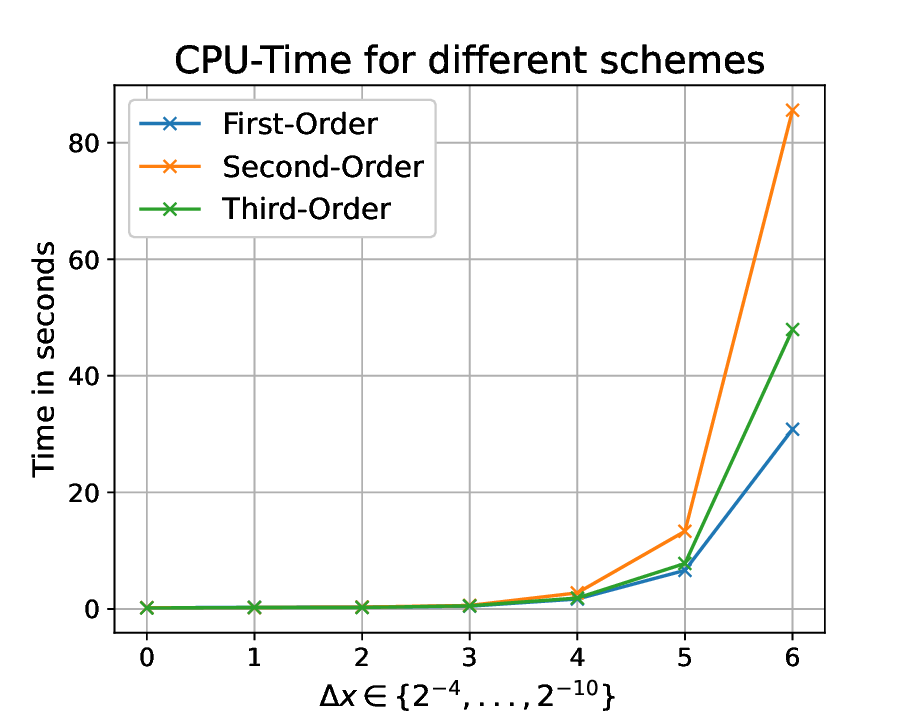}
	\caption{CPU time}
	\label{fig:TwoConsumer_index1_time}
\end{subfigure}
\caption{Error behavior for different spatial discretization and the corresponding CPU time in the index 1 case.}
\label{fig:TwoConsumer_index1}
\end{figure}
It is easy to verify that these parameters and equations satisfy the DAE \eqref{eq:semi_expl} with the underlying differential equations \eqref{eq:endeq}, the coupling conditions \eqref{couplingCond} and the other algebraic equations \eqref{connectNetwork}. 
After some resorting we assume that $z = (z_1, z_2)^\top$, where
\begin{align*}
z_1 &= \begin{pmatrix}
T_1 \\
\vdots \\
T_N
\end{pmatrix}, & z_2 &= \begin{pmatrix}
v_1 \\
\vdots \\
v_N \\
p_{1,0} \\
\vdots \\
p_{N,N}
\end{pmatrix}. 
\end{align*}
The time horizon is given by $[t_0, t_f]$ with $t_0 = 0$ and $t_f=1$ and the discretization is given by
\begin{align*}
0 = t_0 < t_1 < \cdots < t_{n_t} = t_{n_f}=1.
\end{align*}
With $z_{num}$ we denote the numerical solution and with $z_{ana}$ the analytic solution. Then we can compare the error behavior with the discrete absolute, the discrete $L_2$ and the discrete $L_{\infty}$ error with the expressions
\begin{itemize}
\item Discrete $L_1$-error:
\begin{align*}
&\Vert z_{ana} - z_{dis}\Vert_{abs} \\
&\quad= \Delta t \Delta x \sum_{i=0}^{{N_t}} \sum_{j=0}^{{N_x}} \vert z_{1,ana}(t_i,x_j)-z_{1,num}(t_i,x_j) \vert + \Delta t \sum_{i=0}^{{N_t}} \vert z_{2,ana}(t_i)-z_{2,num}(t_i) \vert.
\end{align*} 
\item Discrete $L_2$-error:
\begin{align*}
&\hspace{-5mm}\Vert z_{ana} - z_{dis}\Vert_{L_2}\\
&\hspace{-5mm}= \bigg( 
\Delta t \sum_{i=0}^{{N_t}}\Big(\Delta x\sum_{j=0}^{{N_x}} \vert z_{1,ana}(t_i,x_j)-z_{1,num}(t_i,x_j) \vert^2 + \vert z_{2,ana}(t_i)-z_{2,num}(t_i) \vert^2\Big)
\bigg)^{1/2}.
\end{align*}
\item Discrete $L_{\infty}$-error: 
\begin{align*}
&\hspace{-5mm}\Vert z_{ana} - z_{dis}\Vert_{L_{\infty}}\\
&\hspace{-5mm}= \max_{i = 0, \ldots, N_t} \bigg\{  
\max_{j = 0, \ldots, N_x} \Delta t \Delta x \vert z_{1,ana}(t_i,x_j)-z_{1,num}(t_i,x_j) \vert, 
\Delta t \vert z_{2,ana}(t_i)-z_{1,num}(t_i) \vert
\bigg\}.
\end{align*}
\end{itemize}
\begin{figure}
\begin{subfigure}{.5\linewidth}
	\centering
	\includegraphics[width=\textwidth]{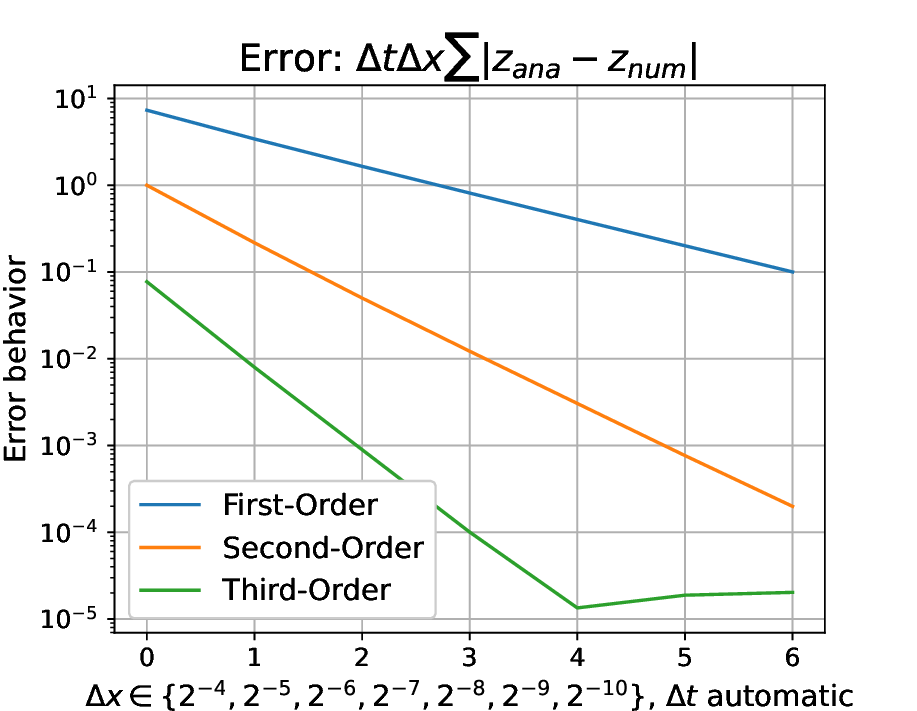}
	\caption{discrete $L_1$-error}
	\label{fig:TwoConsumer_index2_abs}
\end{subfigure}
\begin{subfigure}{.5\linewidth}
	\centering
	\includegraphics[width=\textwidth]{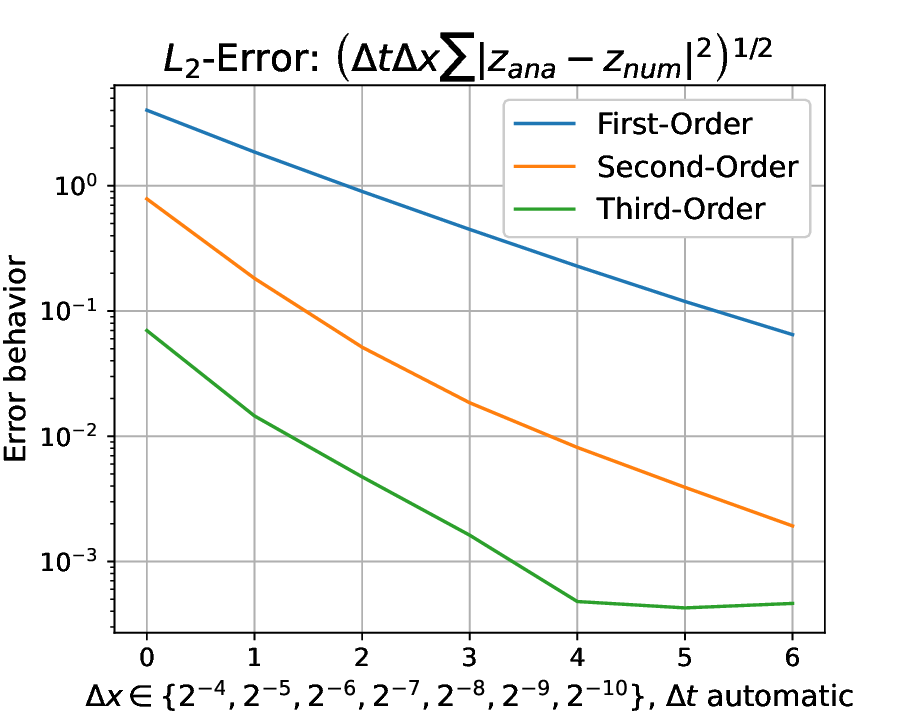}
	\caption{discrete $L_2$-error}
	\label{fig:TwoConsumer_index2_l2}
\end{subfigure}
\begin{subfigure}{.5\linewidth}
	\centering
	\includegraphics[width=\textwidth]{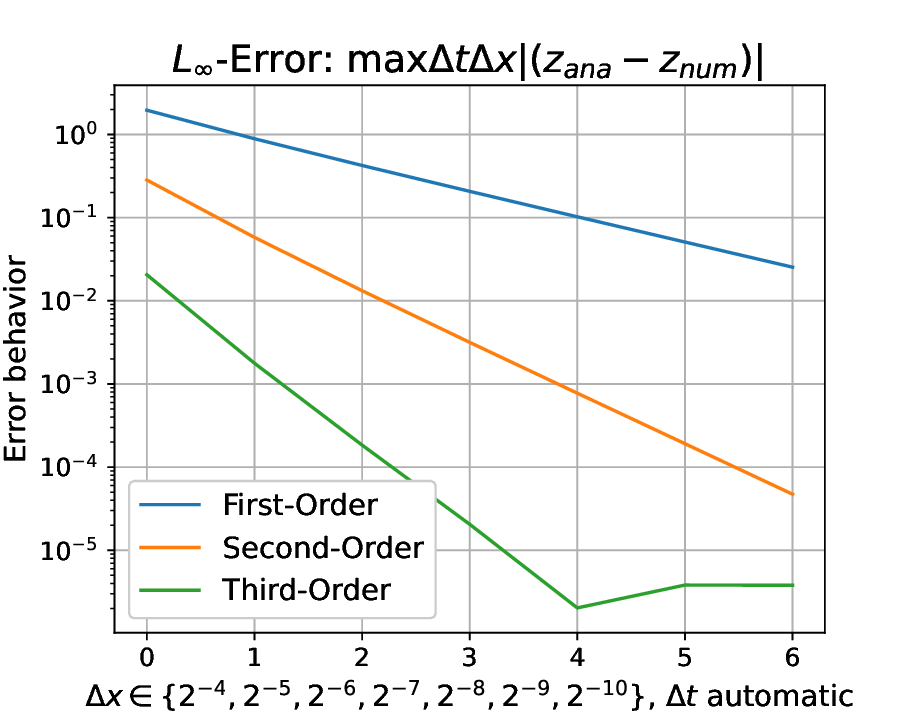}
	\caption{discrete $L_{\infty}$-error}
	\label{fig:TwoConsumer_index2_linf}
\end{subfigure}
\begin{subfigure}{.5\linewidth}
	\centering
	\includegraphics[width=\textwidth]{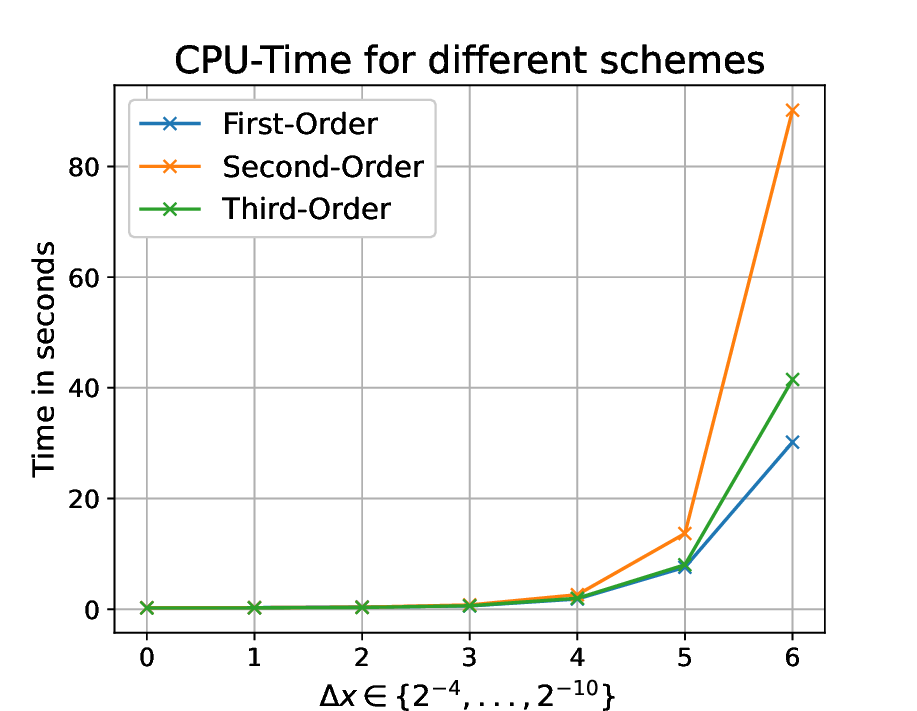}
	\caption{CPU time}
	\label{fig:TwoConsumer_index2_time}
\end{subfigure}
\caption{Error behavior for different spatial discretization and the corresponding CPU time in the index 2 case.}
\label{fig:TwoConsumer_index2}
\end{figure}
In Figure \ref{fig:TwoConsumer_index1} we see the numerical behavior for the reduced model. The expected order is clearly visible in all three cases. The jumps in the second and third order cases are interesting. Here, a numerical artifact seems to have crept in. Also of particular note is the fact that the computation time for the order two case is significantly higher than in the order 3 case. \\
\\
The full model is shown in Figure \ref{fig:TwoConsumer_index2}. Here the expected order is clearly visible in each case, too. It is interesting that the most accurate method (third order) stagnates at an error of about $\approx 10^{-5}$ (or $\approx 10^{-4}$ in the $L_2$ case) and does not improve. The numerical artifact is not seen here. As in the index 1 case, the order 2 method takes more time than the order 3 method although it is less accurate. \\
\\
In summary, neither the index 1 nor the index 2 case is significantly better than the other. The numerical artifact is comparatively insignificant and the computation times do not give each other anything at almost the same accuracy. What can be stated in any case is the fact that it is not worthwhile to simulate with the order 2 method, since the third order method is both more exact and faster (in the sense of CPU time).

\subsection{Simulation in a small network}
\label{sec:4.3}
\begin{figure}
\begin{subfigure}{.7\linewidth}
	\centering
	\includegraphics[width=\textwidth]{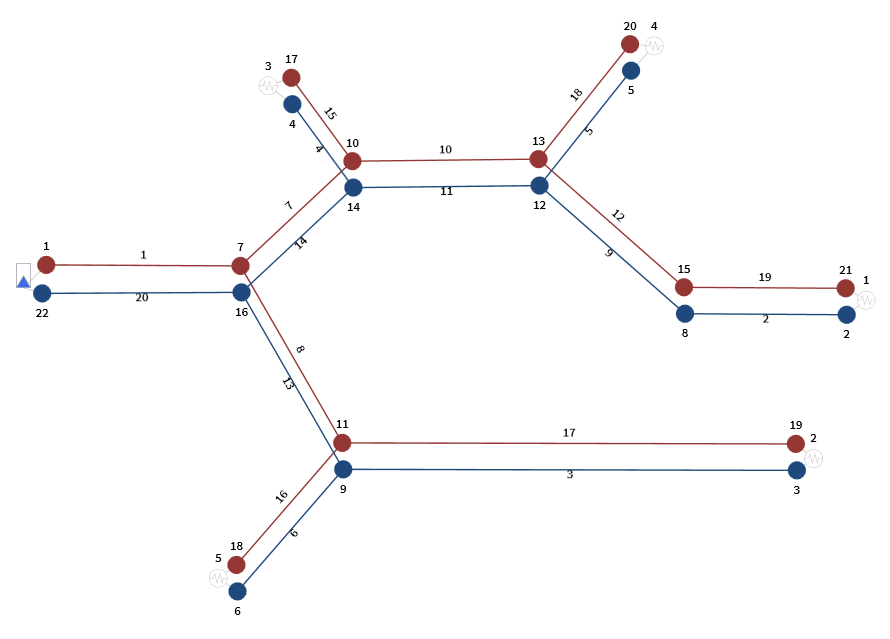}
	\caption{Network with five consumer, $20$ pipes and $22$ nodes.}
	\label{fig:fiveConsumer_network}
\end{subfigure}
\begin{subfigure}{.5\linewidth}
	\centering
	\includegraphics[width=\textwidth]{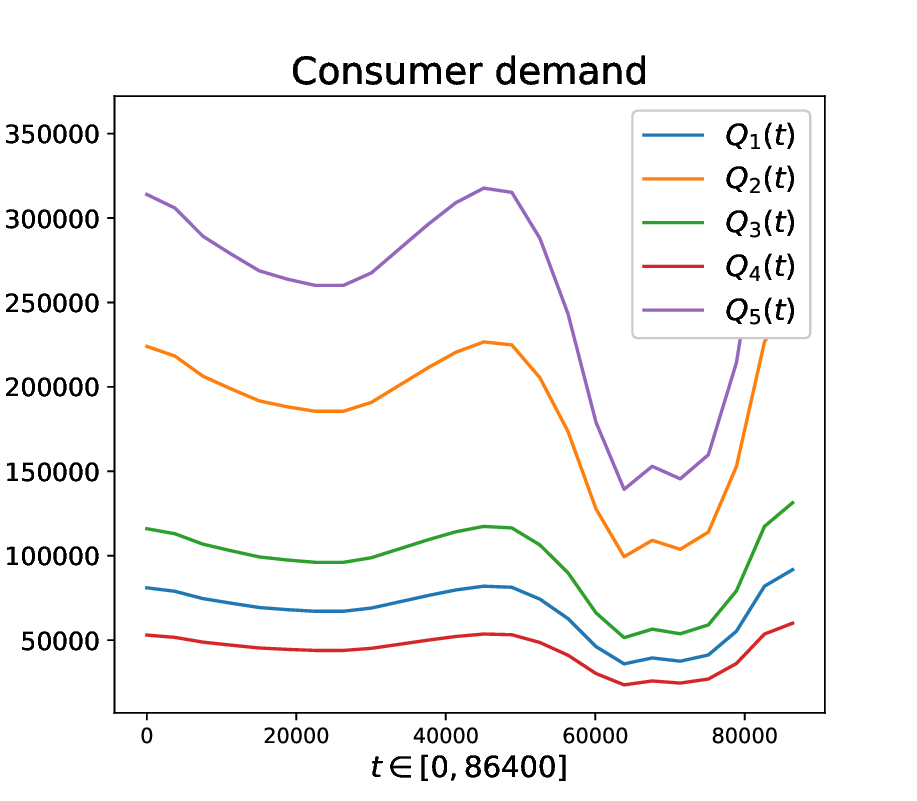}
	\caption{Different consumer demands in the network.}
	\label{fig:fiveConsumer_demand}
\end{subfigure}
\begin{subfigure}{.5\linewidth}
	\centering
	\includegraphics[width=\textwidth]{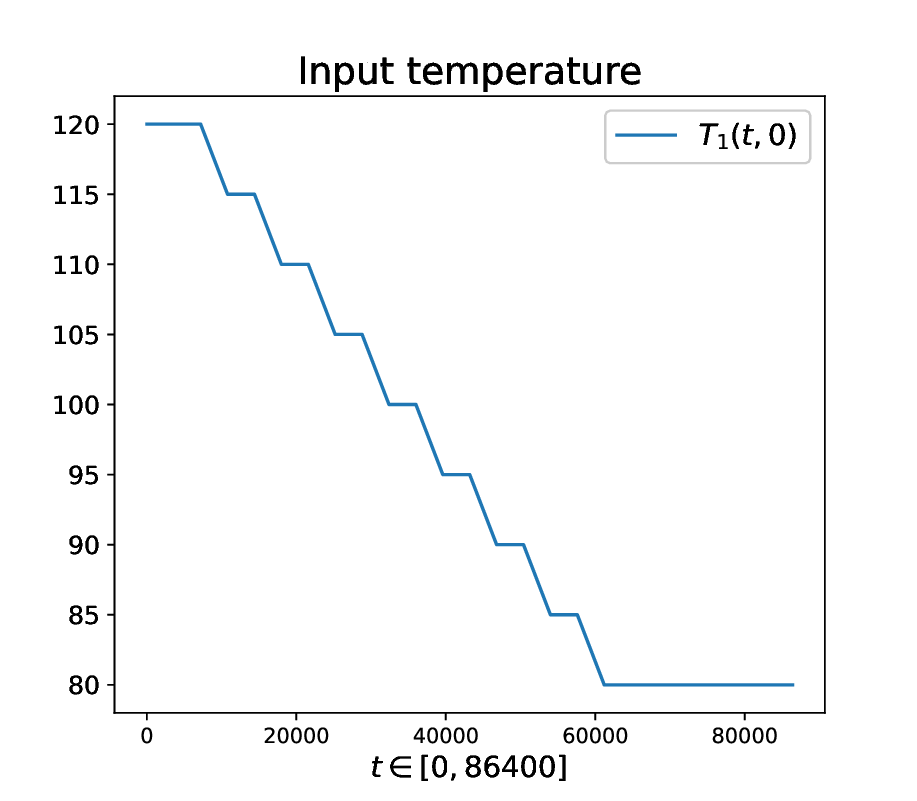}
	\caption{Inlet temperature.}
	\label{fig:fiveConsumer_temperature}
\end{subfigure}
\caption{Some of the input data for the network.}
\end{figure}

In this section, we will use realistic data to simulate a small network. This can be considered as an academic example. The network is shown in Figure \ref{fig:fiveConsumer_network} and we will simulate it with $24$ hours e.g. the time horizon is given by
\begin{align*}
\mathcal{T} = [0,86400].
\end{align*} 
In Figure \ref{fig:fiveConsumer_network} you can see a numbering of the pipes, nodes and consumers. The red lines and the red circles represent the forward flow in the pipes and nodes, the blue ones the return flow back to the depot. We use for all pipes the parameters
\begin{align*}
k &= 0.31, & c_p &= 4160, & d &= 0.1071, & \rho &= 960, \\
T_{ext} &=20, & k_{rough} &= 0.0001, & L &= 100, & \Delta h &= 0.
\end{align*}
The consumption of the consumers can be seen in Figure \ref{fig:fiveConsumer_demand}. This reflects the consumption for one day. Over these $24$ hours we add tempered water to the network as shown in Figure \ref{fig:fiveConsumer_temperature}, starting with $120^\circ\text{C}$ hot water and ending with $80^\circ\text{C}$. In the return flow we assume that the water has always $60^\circ\text{C}$. In addition, we set the inlet pressure to $8$ bar and the pressure arriving at the power plant to $2$ bar. Each pipe is discretized into $100$ segments.
Based on the results in Section \ref{sec:4.1}, we will use the third order discretization here. \\
\\
\begin{figure}
\begin{subfigure}{.5\linewidth}
	\centering
	\includegraphics[width=\textwidth]{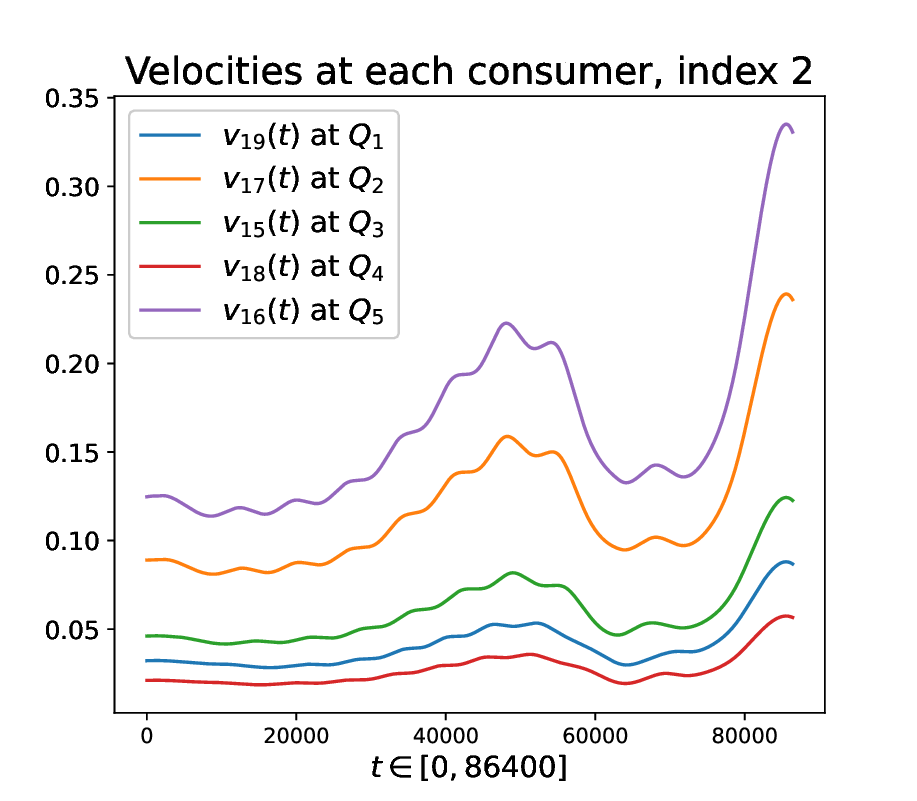}
	\caption{Velocities at the consumers place.}
	\label{fig:fiveConsumer_velo}
\end{subfigure}
\begin{subfigure}{.5\linewidth}
	\centering
	\includegraphics[width=\textwidth]{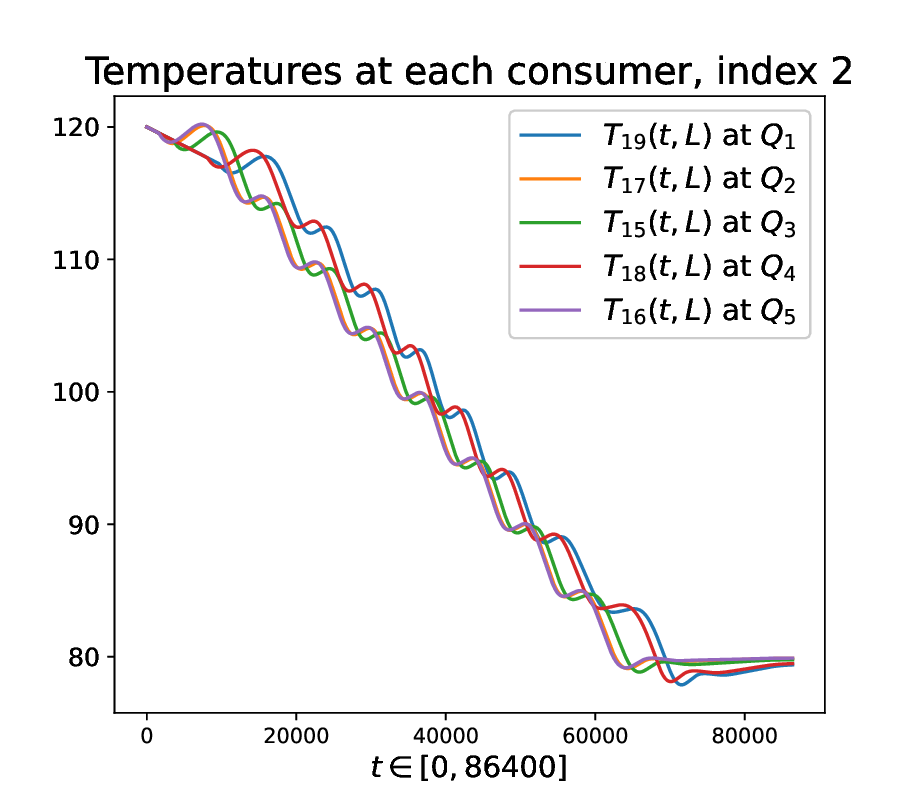}
	\caption{Temperatures at the consumers place.}
	\label{fig:fiveConsumer_temp}
\end{subfigure}
\begin{subfigure}{.7\linewidth}
	\centering
	\includegraphics[width=\textwidth]{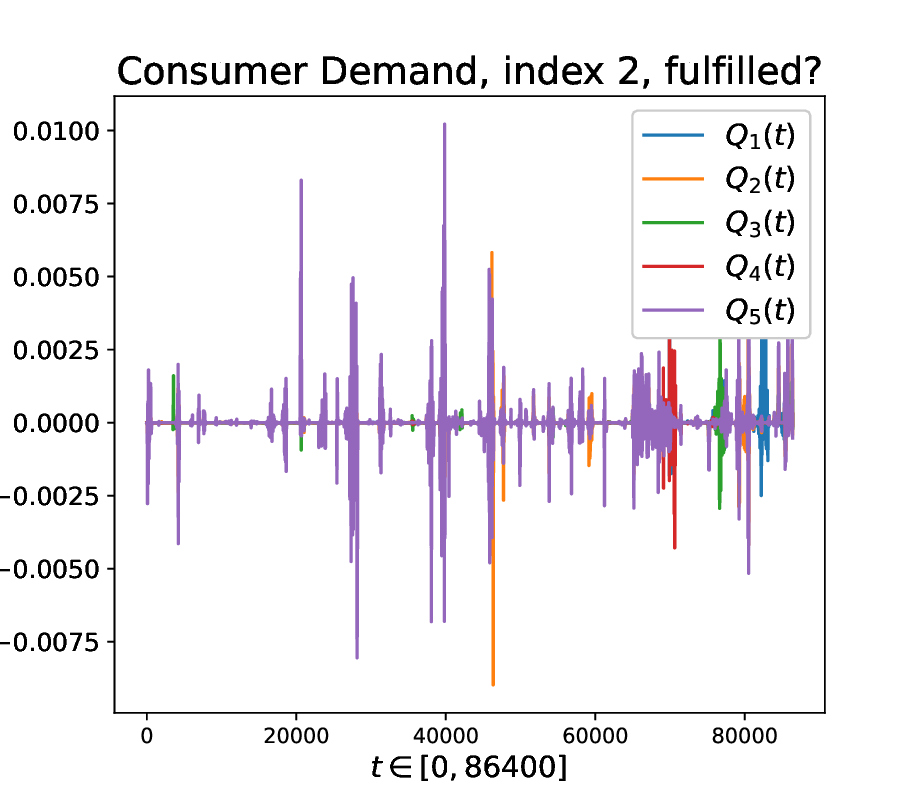}
	\caption{How good is the consumer demand fulfilled.}
	\label{fig:fiveConsumer_cons_fulfilled}
\end{subfigure}
\caption{Some of the simulation results.}
\end{figure}
The differences in the results for the index 1 and index 2 cases are marginal and therefore not discussed further here. Only the CPU time of the index 1 case is a few seconds faster, which is hardly significant in a total time of over two minutes. Therefore, we will now only consider the solutions in the index 2 case. \\
\\
Since the interesting points in a network are essentially the consumers, we take a closer look at them here. In the Figures \ref{fig:fiveConsumer_velo}, \ref{fig:fiveConsumer_temp} and \ref{fig:fiveConsumer_cons_fulfilled} we see with which speed the water arrives at the consumers and with which temperature. Furthermore we see how well the consumer demand equation is fulfilled.
A particularly interesting observation is that the velocities essentially follow the consumption, compare Figure \ref{fig:fiveConsumer_demand} and Figure \ref{fig:fiveConsumer_velo}. But with the velocities one see additionally small bumps which show additionally the sinking temperature and thus the smaller becoming energy. With higher velocities, this loss of energy must be compensated for so that the consumers are sufficiently supplied.
We can make similar statements about the temperature in the network. With a slight delay, essentially the same temperature arrives at the consumers that was previously pumped into the network at the depot. We can also see slight dents here, which is also due to the changing temperature.
All in all, it can be stated that in this small example, with data that is realistic for this size, we get results from the simulation process which are expected.

\subsection{Simulation part of a real network}
\label{sec:4.4}
\begin{figure}[h]
\centering
\includegraphics[width=\textwidth]{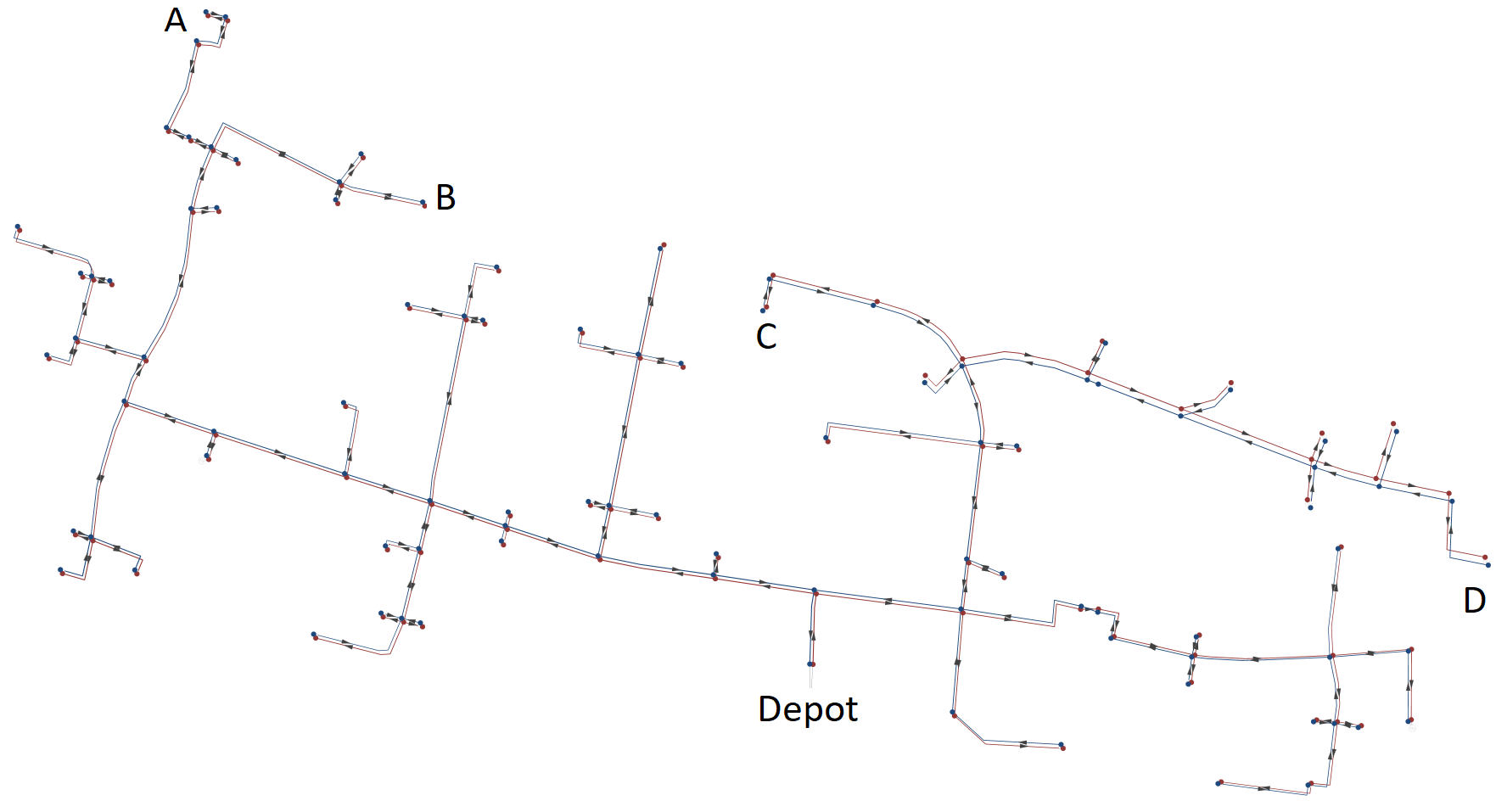}
\caption{Part of a real network.}
\label{fig:big_net}      
\end{figure}

\begin{figure}
\begin{subfigure}{.33\linewidth}
	\centering
	\includegraphics[width=\textwidth]{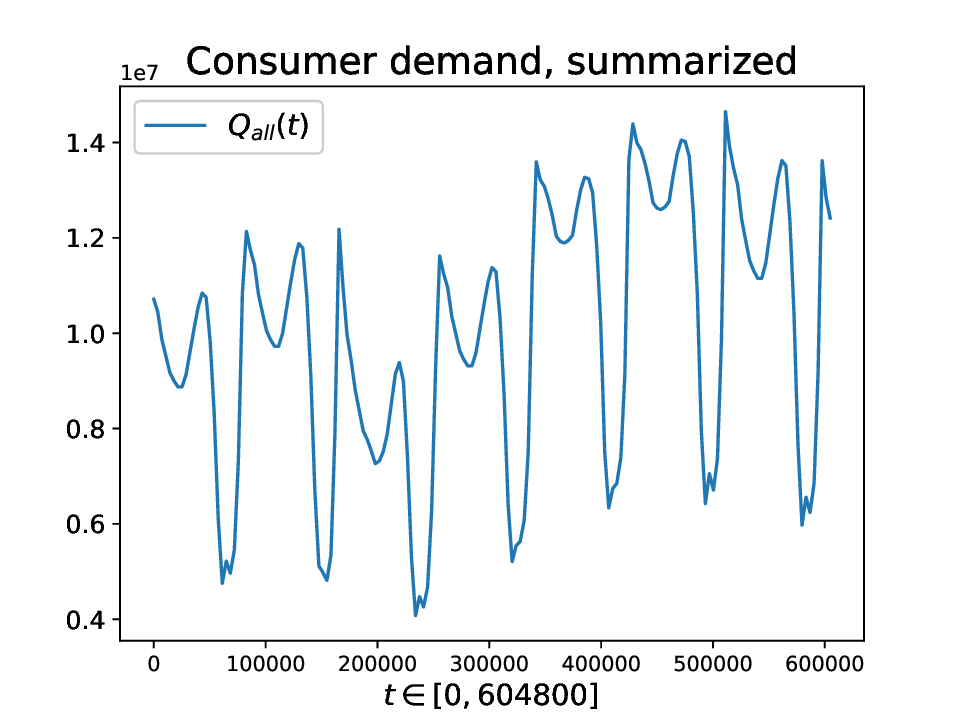}
	\caption{Total consumption of all consumers for one week.}
	\label{fig:total_consumption}
\end{subfigure}
\begin{subfigure}{.33\linewidth}
	\centering
	\includegraphics[width=\textwidth]{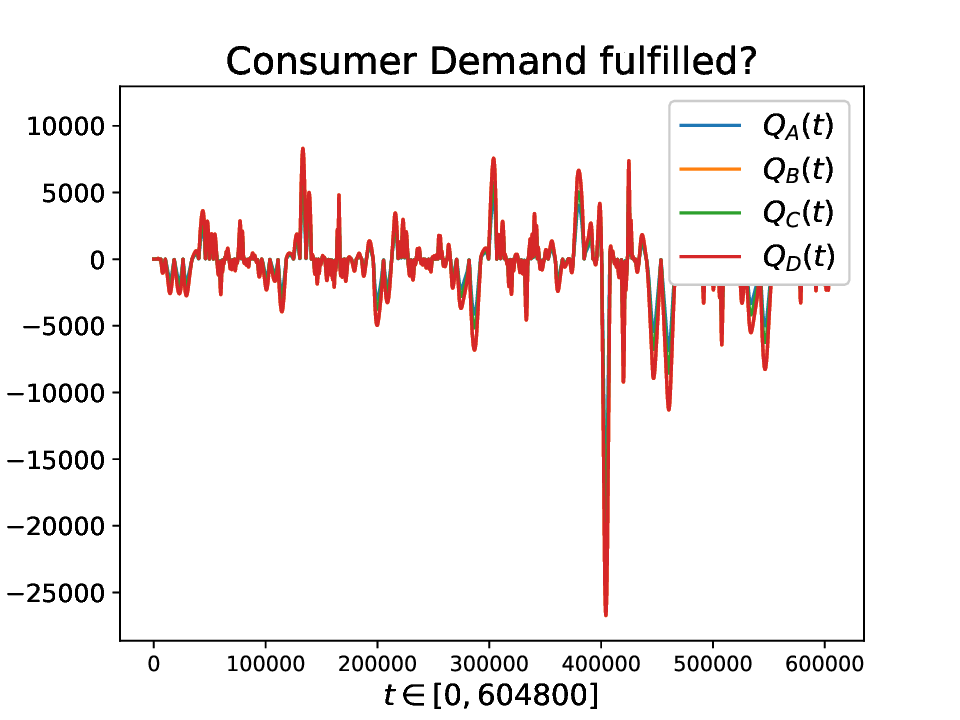}
	\caption{How good is the consumer demand fulfilled, index 1.}
	\label{fig:big_net_con_ful_index1}
\end{subfigure}
\begin{subfigure}{.33\linewidth}
	\centering
	\includegraphics[width=\textwidth]{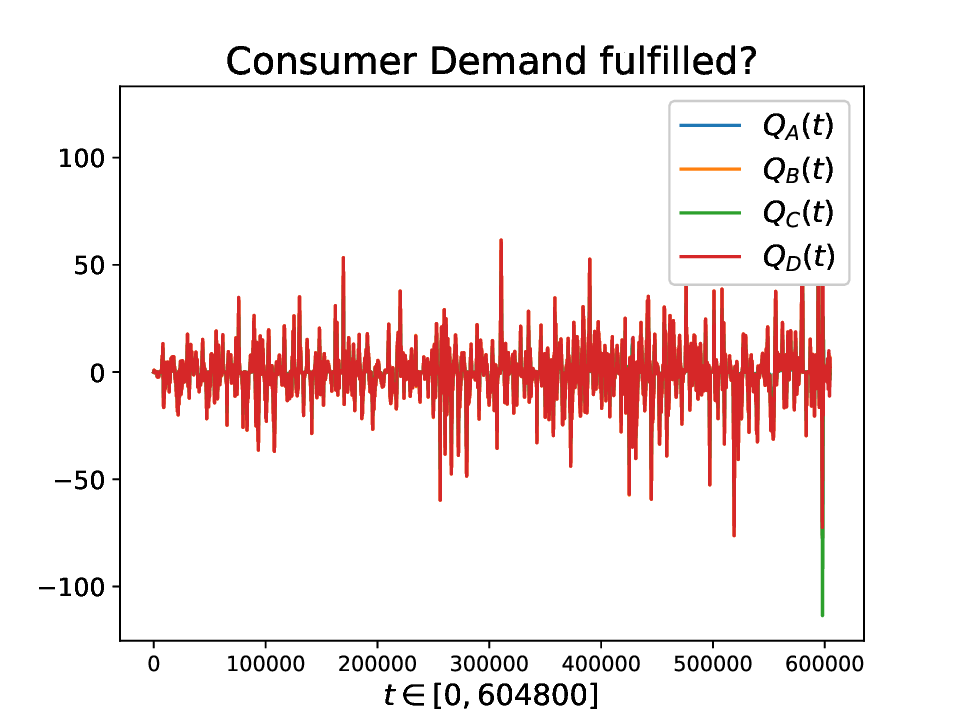}
	\caption{How good is the consumer demand fulfilled, index 2.}
	\label{fig:big_net_con_ful_index2}
\end{subfigure}
\end{figure}

In this last example, we simulate a part of a real network. The data is provided by our project partners (Rechenzentrum für Versorgungsnetze Wehr GmbH\footnote{See \url{https://www.rzvn.de}}). The layout of the network can be seen in Figure \ref{fig:big_net}.
The network consists of a total of $193$ nodes and $51$ consumers. Of the $190$ pipes, with a total length of $7988$ meters, $95$ each of the pipes are inflow and return. The total consumption of all consumers for one week can be seen in Figure \ref{fig:total_consumption}. The following operating case is simulated. We assume that the depot pumps $100^\circ$C hot water with $6.8$ bar into the network. The return flow is $70^\circ$C and arrives at the depot with $1.2$ bar. Since the pipes are buried in the ground, we assume a constant ambient temperature of $6^\circ$C. 
\\
\\
We use a coarse space discretization and the pipes are discretized differently. If a pipe is less than $40$ meters long, it is discretized into $5$ sections and all pipes longer than $40$ meters are discretized into $20$ sections.
\\
\\
As a result, the IDA solver must solve a nonlinear system with almost $3000$ unknowns at each time step. The solver uses an automatic step size control and we use as relative tolerance $10^{-4}$ and as absolute tolerance $10^{-6}$. We compare some results for the index $1$ and index $2$ case. For this we use the same initial value. Note that a consistent initial value for the index $2$ case is also consistent for the index $1$ case. Therefore, we calculate it as described above for the index $2$ case and use it for both cases. This calculation takes about $13$ seconds.
For the simulation in the reduced index $1$ case the solver needs about $28$ seconds, for the index $2$ case about $69$ seconds.  

An interesting point in such a network are consumers far away from the depot. Here we consider the customers A, B, C and D see Figure \ref{fig:big_net}. How well the consumption equation is fulfilled in each case can be seen in Figures \ref{fig:big_net_con_ful_index1} and \ref{fig:big_net_con_ful_index2}. Here a clear difference between the index $1$ and index $2$ case can be seen. Although we use the same initial value and the same tolerances for the solver, in the index $1$ case the equations are clearly worse fulfilled.
\\
\begin{figure}
\begin{subfigure}{.33\linewidth}
	\centering
	\includegraphics[width=\textwidth]{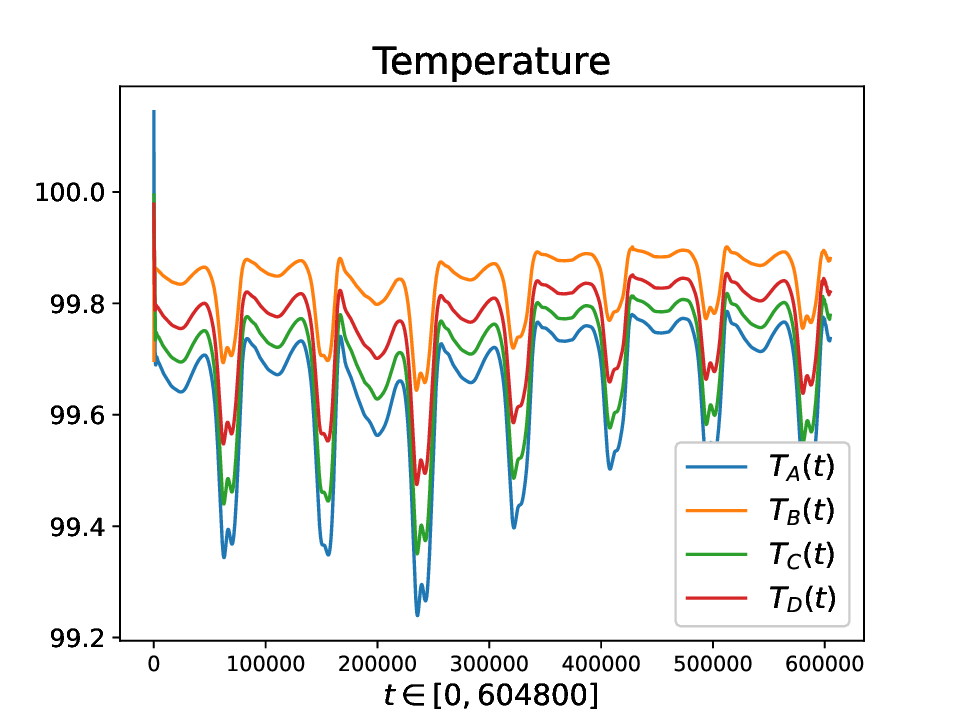}
	\caption{Different temperatures at consumers A, B, C and D.}
	\label{fig:big_net_temperature}
\end{subfigure}
\begin{subfigure}{.33\linewidth}
	\centering
	\includegraphics[width=\textwidth]{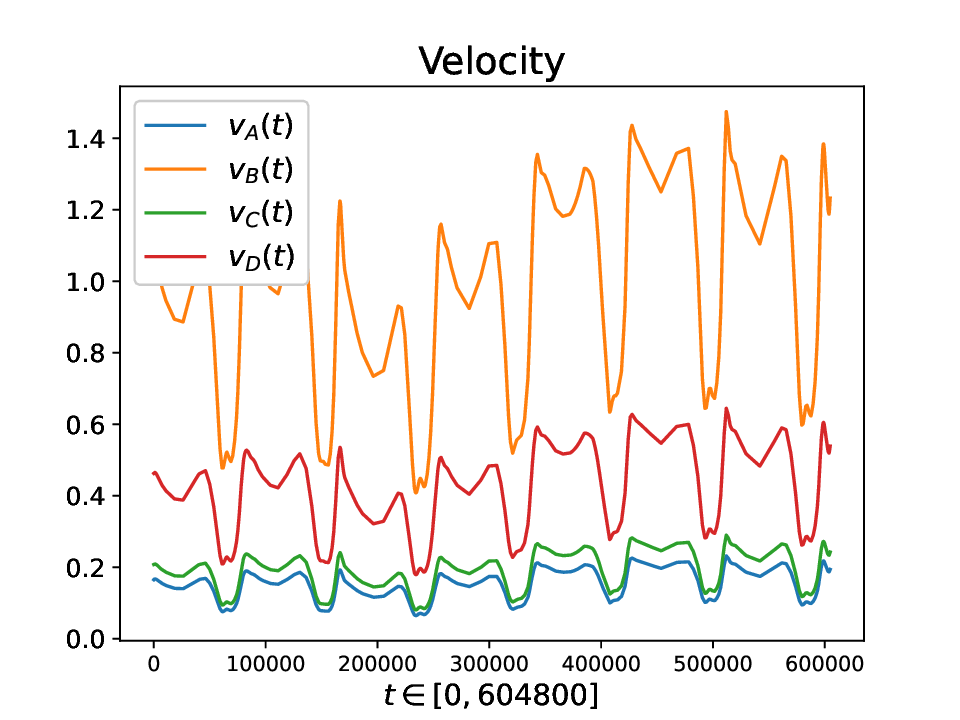}
	\caption{Different velocities at consumers A, B, C and D.}
	\label{fig:big_net_velocity}
\end{subfigure}
\begin{subfigure}{.33\linewidth}
	\centering
	\includegraphics[width=\textwidth]{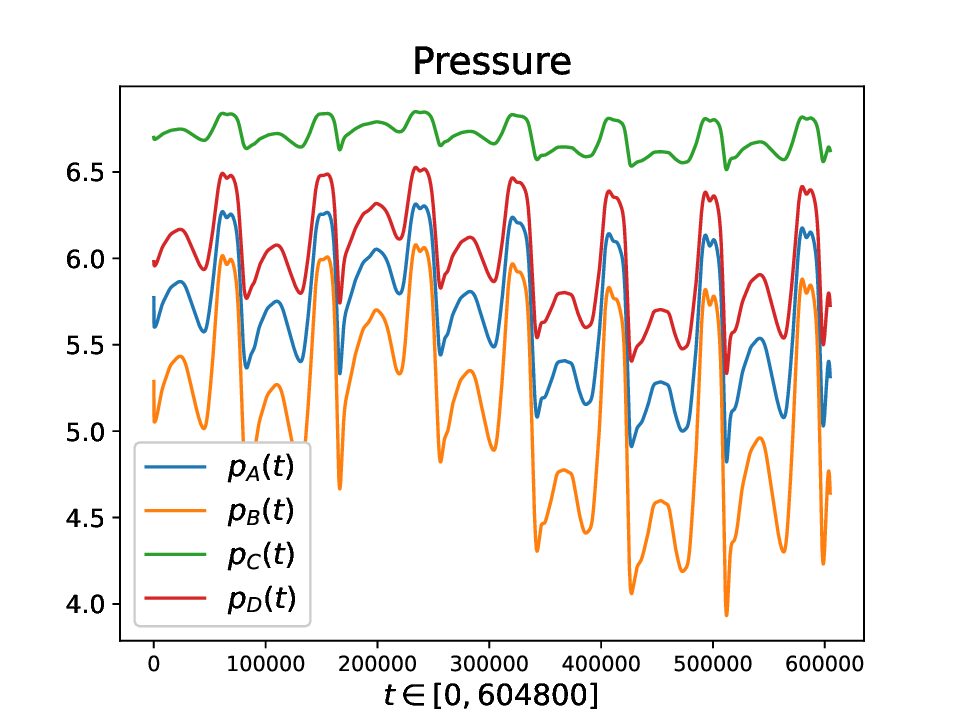}
	\caption{Different pressures at consumers A, B, C and D.}
	\label{fig:big_net_pressure}
\end{subfigure}
\end{figure}
\\
Figure \ref{fig:big_net_temperature} shows the temperatures that arrive at the respective consumers. These are close to the $100^\circ$C that is pumped into the grid. It is easy to see that there is hardly any heat loss. The different deflections caused by the consumers can also be seen very clearly. In the times in which the consumers consume little, the temperature also goes down somewhat and vice versa. This is related to the flow velocity. If the consumers do not pull much, the velocity is reduced; this can also be seen clearly, compare Figure \ref{fig:big_net_velocity}. Finally, the velocities result from the pressure differences in the respective pipes. The corresponding pressure can be seen in Figure \ref{fig:big_net_pressure}, which also shows the same qualitative behavior as the consumers.
\\
\\
Finally, we compare some values in the Euclidean norm below. If we denote with $z_{full}$ the solution of the index $2$ system and with $z_{red}$ the solution of the index $1$ problem we get as difference
\begin{align*}
\sqrt{\sum_{k=1}^{n_t} \Delta t \Vert z_{full}^k - z_{red}^k \Vert_2^2} = 0.0908
\end{align*}
As a reminder, the only difference between the full and the reduced model is that we neglect the term $\dot{v}$ in the modified incompressible euler equation \eqref{eq:euler}. Since the velocity results essentially from the pressure difference, let us consider the difference in the Euclidean norm of velocities and pressures at the consumers A, B, C and D:
\begin{center}
\begin{tabular}{c|c|c|c|c}
 & A & B & C & D \\ \hline 
Velocity & $0.13$ & $0.87$ & $0.17$ & $0.38$ \\
Pressure & $1.35$ & $1.94$ & $0.31$ & $1.08$ \\
\end{tabular}
\end{center}
If we compare the difference of velocities and pressures at the consumers A, B, C and D with a relative discrete $L_2$ Norm, e.g.
\begin{align*}
\frac{\Bigl( \int_{t_0}^{t_f} \bigl( h_{full}(t) - h_{red}(t)\bigr)^2 dt \Bigr)^{1/2} }{\Bigl( \int_{t_0}^{t_f} \bigl( h_{full}(t)\bigr)^2 dt \Bigr)^{1/2}}
\end{align*} 
we get:
\begin{center}
\begin{tabular}{c|c|c|c|c}
 & A & B & C & D \\ \hline 
Velocity & $1.6e-3$ & $1.6e-3$ & $1.6e-3$ & $1.6e-3$ \\
Pressure & $4.5e-4$ & $7.2e-4$ & $87.3e-4$ & $3.4e-4$ \\
\end{tabular}
\end{center}
Note that we used a trapezoidal rule for the integral. These numbers, especially the relative error, show that the term $\dot{v}$ in the modified incompressible euler equation \eqref{eq:euler} is negligible. In fact, the error is smaller or the same as the stepsize for the time integration. Together with the fact that for small networks it makes hardly any difference in terms of CPU time whether the full model or the reduced model is solved and for large networks the reduced model is solved significantly faster, it can be stated here that it makes perfect sense to solve the reduced model. However, if you need precise results for the consumer equations, it may be advisable to accept the possibly higher computing time and solve the full index $2$ model.

\section{Conclusions}
\label{sec:5}

In this article we could show that our model works well and the results are convincing. In the case that the district heating network is represented as a graph-theoretic tree, we could show existence and uniqueness of the solution. \\
\\
In numerics, we looked at different discretizations of space and compared them for a constructed example where we know the exact solution. We were able to show that more accurate results can be expected for higher order, but that these require more computing time. Surprisingly, the second order discretization method requires more CPU time than the third order one.\\
\\
Our model is essentially based on the work of \cite{Borsche19} and \cite{Krug19}. However, this is an index $2$ DAE which is numerically more complex to handle. Instead of going the classical way of index reduction and encountering the possible problems, we followed \cite{rein2019model} and neglected the term $\dot{v}$ which leads to an index $1$ DAE. In the simulation we have seen that for practical examples the differences are marginal and negligible. Finally, we have seen that the simulation is also feasible for parts of real networks.

\section*{Appendix A. Functions $f_1, \, f_2, \ g_1$ and $g_2$ modeling}
The main objective of this section is the modeling of the functions $f_1$, $f_2$, $g_1$ and $g_2$ in a abstract way to handle any network automatically in a consistent manner. For this we assume, that we have a
directed acyclic graph $\g = (\n,\e)$ with direction following ordering. Let $\mathcal{I}$ denote its incidence matrix and $\c^1, \, \c^2$ its consumer matrix. Let $N$ be the number of pipes, $\oN$ the 
number of nodes with $\oN = n_s+n_d+n_{junc}$, $n_i $ the number of discretization points in pipe $e_i$ for $i = 1,...,N$ $n$ the sum of all $n_i$ with $n := \sum\limits_{i = 1}^N n_i$. We denote $p_{N,L}$ the given pressure for the demand node for the power plant
and $T_2^A,...,T_{n_s}^A$ the given temperature for the supply nodes before a consumer. Moreover $Q_1,...,Q_{n_c}$ describe the consumers. For a better reading we left out the arguments for $t$. 
Let $\mathcal{I}_+ := \max \{\mathcal{I},0\}$ and $\mathcal{I}_- := \min \{\mathcal{I},0\}$ be the positive and negative part of the incidence matrix $\mathcal{I}$ and $\mathcal{I}^r := (\mathcal{I}_{i,j})_{m_s+1 \leq i \leq m_s+m_0, 1\leq j \leq N}$ the 
reduced incidence matrix. \\
\\
In order to write our system in the form \eqref{eq:semi_expl}  we do a spatial discretization (forward Euler) by the method of lines for the transport equation \eqref{eq:transport} and get 
\begin{align}
\label{discreteTemp}
\begin{pmatrix}
\dot{T}_{i,1} \\
\vdots \\
\dot{T}_{i,n} \\
\end{pmatrix} = - \frac{v_i}{\Delta x_i} \begin{pmatrix}
T_{i,1} \\
T_{i,2} - T_{i,1} \\
\vdots \\
T_{i,n_i} - T_{i,n_i-1} 
\end{pmatrix} - \frac{4 k}{c_p d_i \rho} \left[ \begin{pmatrix}
T_{i,1} \\
\vdots \\
T_{i,n_i}
\end{pmatrix} - \begin{pmatrix}
T_{ext} \\
\vdots \\
T_{ext}
\end{pmatrix} \right],
\end{align}
for $i = 1,...,N$. Similar we get for the Euler equation \eqref{eq:euler}
\begin{align}
\label{euler2}
\dot{v}_i = -\left(\frac{p_i(L)-p_i(0)}{L_i \bar{\rho}}+\frac{\lambda_i}{2d_i} v_i^2 + g(\partial_x h)\right),
\end{align}
for $i = 1,...,N$.
For each pipe $e_i \in \e$ we have the set of variables $T_{i,1},...,T_{i,n_{i}}$, $v_i$ and $p_i(L),p_i(0)$ depending on $t$ that represent the discrete analog of $T,v,p$, respectively. 
Here $n_i$ is the number of discretization points in pipe $e_i$.
The coupling conditions \eqref{couplingCond} leads to the following equations:
\begin{subequations}
\label{discouplingCond}
\begin{align}
\label{disconserveMass}
\sum\limits_{e_i \in \sigma_j\cup \Sigma_j} A_i \rho v_i = 0, \\
\label{disconserveEnergy}
\sum\limits_{e_i \in \sigma_j}c_p A_i \rho v_i T_{i,n_i}- \sum\limits_{e_i \in \Sigma_j}c_p A_i \rho v_i T_{i,1}= 0, \\
\label{discontpressure}
p_i(L) = p_l(0) \mbox{ for all } i \in \sigma_j, \, l \in \Sigma_j, \\
\label{disperfectMix}
T_{i,1} = T_{l,1} \mbox{ for all } i,l \in \Sigma_j,
\end{align}
\end{subequations}
where $\sigma_j$ is the set of pipes incoming the node $j$ and $\Sigma_j$ the set of pipes outgoing of node $j$.
For consumer $k$ with incoming pipe $e_i$ and outgoing pipe $e_l$ it holds
\begin{subequations}
\label{disconnectNetwork}
\begin{align}
\label{disnoMassLost}
q_i &= q_l, \\
\label{disconsumerDemand}
Q_k(t) &= c_p A_i q_i(T_{i,n_i} - T_{out}).
\end{align}
\end{subequations}
Notice, that $T_{out}$ is a fixed number.

\paragraph{Functions $f_1$ and $f_2$}
\mbox{}\\
We start with the functions $f_1$ and $f_2$ which describes the discretized PDE and the ODE. For the PDE \eqref{discreteTemp} we define the matrices 
\begin{align*}
A_{f_1}& = \left(\begin{array}{c c c}
-\begin{pmatrix}
-\frac{1}{\Delta x_1} \\ \vdots \\ -\frac{1}{\Delta x_1}
\end{pmatrix}_{\tilde{n}_1} & &  \\
& \ddots &  \\
& & -\begin{pmatrix}
-\frac{1}{\Delta x_N} \\ \vdots \\ -\frac{1}{\Delta x_N}
\end{pmatrix}_{\tilde{n}_N} 
\end{array} \right),\\
\end{align*}
\begin{align*}
A_{f_2} &= \begin{pmatrix}
A_{\tilde{n}_1} & & \\
 & \ddots & \\
 & & A_{\tilde{n}_N}
\end{pmatrix} \text{ with } A_{\tilde{n}_i} = \begin{pmatrix}
1 & & &  \\
-1 & 1 & & \\
 & \ddots & \ddots & \\
 &  & - 1  & 1
\end{pmatrix} \in \R^{\tilde{n}\times \tilde{n}}
\end{align*}
\begin{align*}
 A_{f_4} &= \begin{pmatrix}
- \frac{4k}{d_1 c_p \rho} \eins_{\tilde{n}_1}  & & \\
& \ddots  & \\
& &  - \frac{4k}{d_N c_p \rho} \eins_{\tilde{n}_N}
\end{pmatrix},\\
D_{f_1}& = \left(\begin{array}{c c c | c}
-\begin{pmatrix}
1 \\ \vdots \\ 0
\end{pmatrix}_{\tilde{n}_1} & & & \\
& \ddots & & 0_{\tilde{n} \times 2N} \\
& & -\begin{pmatrix}
1 \\ \vdots \\ 0
\end{pmatrix}_{\tilde{n}_N} &
\end{array} \right),
\text{ and } d_1 = \begin{pmatrix}
\frac{4k T_{ext}}{d_1 c_p \rho} \\
\vdots \\
\frac{4k T_{ext}}{d_1 c_p \rho} \\
\vdots \\
\frac{4k T_{ext}}{d_N c_p \rho} \\
\vdots \\
\frac{4k T_{ext}}{d_N c_p \rho}
\end{pmatrix}.
\end{align*}
This leads to $\dot{x_1}(t) = f_1(t, x(t), y(t))$ with 
\begin{align*}
f_1(t, x(t), y(t)) = A_{f_1} x_2(t)  .* \left( A_{f_2}x_1(t) + D_{f_1}y(t) \right) +  A_{f_4} x_1(t) + d_1
\end{align*}
and
\begin{align*}
A_{f_1}  \in \R^{\tilde{n} \times N},  \,
A_{f_2}, A_{f_4}  \in \R^{\tilde{n} \times \tilde{n}},  \,
D_{f_1}  \in \R^{\tilde{n} \times 3N}, \,
d_1  \in \R^{\tilde{n} }.
\end{align*}
\begin{remark}
\label{rem:spatial_discretization}
In some cases we want to use another scheme for the spatial discretization. In particular we have used the following three different discretizations.
\begin{align*}
\partial_x T(t,x+h) = 
\begin{cases}
\frac{1}{h}\left( T(t,x+h)-T(t,x) \right) , \quad & \text{order=1}, \\
\frac{1}{2h}\left( T(t,x+2h)-T(t,x) \right) , \quad & \text{order=2}, \\
\frac{1}{6h}\left( T(t,x-h)-6T(t,x) + 3T(t,x+h)+2T(t,x+2h) \right) ,& \text{order=3}. \\
\end{cases}
\end{align*}

To get the corresponding points on the boundary we need some ghost points which should be calculated with the same order we use here. For more details about ghost points at boundary see e.g. \cite{Albaiz14}. The idea is to interpolate with the appropriate order e.g \begin{align*}
T(x) & \approx a + bx & & \text{ for order } = 2, \\
T(x) & \approx a + bx + cx^2 & & \text{ for order } = 3.
\end{align*}

Assume that we want to discretize T with $N$ spatial discretization points. In the case that we use a second order scheme we need therefore a ghost point for $T_{N+1}$. Simple calculations and insert into the corresponding equations shows that we need to change only $A_{\tilde{n}_i}$, e.g.
\begin{align*}
\begin{pmatrix}
1 & & & & \\
-1 & 0 & 1 & & \\
& & \ddots & \ddots & \\
& & & -2 & 2
\end{pmatrix}
\end{align*}
In the case that we use a third order scheme we need to solve for the ghost point on the left boundary.
\begin{align*}
T_0 &= a,  & & a = T_0, \\
T_1 &= a + b \Delta x + c \Delta x^2,  & & b = \frac{-3 T_0+ 4 T_1 - T_2}{2 \Delta x},  \\
T_1 &= a + b 2\Delta x + c 4\Delta x^2, & & c =  \frac{T_0 - 2 T_1 + T_2}{2 \Delta x^2}.
\end{align*}
Therefore our ghost point is given by
\begin{align*}
T_{-1} = 3 T_0 - 3T_1 + T_2.
\end{align*}
After inserting into the correct equation we get for the first two equations
\begin{align*}
dT_1 &= ... \frac{1}{6h}(-3T_0 + 3T2), \\
dT_2 &= ... \frac{1}{6h}(T_0 - 6T_1 + 3T_2 + 2T_3).
\end{align*}
Similar calculations show that we have for the last equation:
\begin{align*}
dT_N &= ... \frac{1}{6h}(3 T_{N-2} - 12T_{N-1} + 9T_{N})
\end{align*}
Therefore we need to change $A_{\tilde{n}_i}$ and $D_{f_1}$ to
\begin{align*}
A_{\tilde{n}_i} = \begin{pmatrix}
0 & 3 & & & & \\
1 & -6 & 3 & 2 & & \\
 & \ddots & \ddots & \ddots & \ddots & \\
 & & 1 & -6 & 3 & 2 \\
 & & & 3 & -12 & 9
\end{pmatrix}, \qquad 
D_{f_1} = \left(\begin{array}{c c c | c}
-\begin{pmatrix}
-3 \\ 1 \\ \vdots \\ 0
\end{pmatrix}_{\tilde{n}_1} & & & \\
& \ddots & & 0_{\tilde{n} \times 2N} \\
& & -\begin{pmatrix}
-3 \\ 1 \\ \vdots \\ 0
\end{pmatrix}_{\tilde{n}_N} &
\end{array} \right).
\end{align*}
\end{remark}
After we finished for the PDE we now define the matrices for the ODE \eqref{euler2} as follows:
\begin{align*}
A_{f_3} = \begin{pmatrix}
 -\frac{\lambda_1}{2 d_1} & & \\
 & \ddots & \\
 & & -\frac{\lambda_N}{2 d_N}
\end{pmatrix}, \quad D_{f_2} = \left( \begin{array}{c|ccccc}
& \frac{1}{L_1 \rho} & -\frac{1}{L_1 \rho} & & &  \\
0_{N \times N} & & & \ddots & & \\
& & & & \frac{1}{L_N \rho} & -\frac{1}{L_N \rho}
\end{array} \right) \text{ and } d_2 = \begin{pmatrix}
- g(\partial x h) \\
\vdots \\
- g(\partial x h) 
\end{pmatrix}.
\end{align*}
This leads to $\dot{x_2}(t) = f_2(t, x(t), y(t))$ with 
\begin{align*}
f_2(t, x(t), y(t)) = A_{f_3} x_2(t)  .* x_2(t)  +  D_{f_2} y(t) + d_2
\end{align*}
and
\begin{align*}
A_{f_3} \in \R^{N \times N}, \,
D_{f_2} \in \R^{N \times 3N}, \,
d_2 \in \R^{N }.
\end{align*}
Using $x = \begin{pmatrix}
x_1 \\ x_2
\end{pmatrix}$ we get
\begin{align*}
\dot{x}(t) = f(t, x(t), y(t)) = \begin{pmatrix}
f_1(t, x(t), y(t)) \\
f_2(t, x(t), y(t))
\end{pmatrix}
\end{align*}

\paragraph{The function $g_1$}
\mbox{}\\
In this section we want describe who we model the function $g_1$.
We start with the equations for the continuity of the pressure and the perfect mixing \eqref{discontpressure} and \eqref{disperfectMix} where we need the matrices $D_{g_1,p}$ and $D_{g_1,T}$. They are given by Algorithm \ref{algo:D_g1P} and \ref{algo:D_g1T}. Notice, that $n_p$ and $n_{out}$ are the number of rows of $D_{g_1,p}$ and $D_{g_1,T}$. After generating them we can set $n_{g_1} = n_p + n_{out} + 2 + n_c + 1$ and 
\begin{align*}
D_{g_1,p} = \begin{pmatrix}
0_{n_p \times N} & D_{g_1,p} \\ \hline
 0_{(n_{g1}-n_p) \times 3N} &
\end{pmatrix} \text{ and } D_{g_1,T} = \begin{pmatrix}
0_{n_p \times 3N} & \\
D_{g_1,T} & 0_{n_{out} \times 2N} \\
0_{(n_{g1}-n_{out}-n_p) \times 3N} &
\end{pmatrix}
\end{align*}
\begin{algorithm}
	\SetAlgoLined
    	\SetKwInOut{Output}{Return}
		\Output{matrix $D_{g_1,p}$}
	Set $a \in \R^\oN$ with $a_i := \sum\limits_{j = 1}^N (\mathcal{I}_+ - \mathcal{I}_-)_{i,j}$ \\
	Set $D_{g_1,p}=[]$ \\
	\For{$i =0,1,...,n_{junc}$}{
		Set $i := n_s + i$;\\
		Set ind $:= nonzero(\mathcal{I}[i,:])$, the nonzero elements in $\mathcal{I}$. \\
		Set $l := a_i$ \\
		\For{$j = 1,...,l-1$}{
			Set tmp = zeros($2m$) \\
			Set s, q = ind[0][j], ind[0][j+1] \\
			\eIf{$\mathcal{I}[i, s] == -1$}{
			$tmp[2*s+1]=-1$}{
			$tmp[2*s]=-1$ }
			\eIf{$\mathcal{I}[i, q] == -1$}{
			$tmp[2*q+1]=1$}{
			$tmp[2*q]=1$} 
			Set $D_{g_1,p}.append(tmp)$ }}
		\caption{Generate $D_{g_1,p}$ for $g_1$}
	\label{algo:D_g1P}
\end{algorithm}
\begin{algorithm}
	\SetAlgoLined
    	\SetKwInOut{Output}{Return}
		\Output{matrix $D_{g_1,T}$}
	Set $\mathcal{I}_{tmp} = \mathcal{I}_{r_p}[count nonzero(\mathcal{I}_r == 1, axis=1) > 1]$
	Set $n_{temp}, m_{temp}=shape(\mathcal{I}_{temp})$
	Set $D_{g_1,T}=[]$ \\
	\For{$i = 0,1, ..., n_{temp}$}{
			 $temp = where(\mathcal{I}[i, :] == 1)$ \\
			 \For{$j = 0,1,...,len(temp)-1$}{
			 	$temp2=zeros(m_{temp}$ \\
			 	$temp2[temp[0]] = 1$ \\
			 	$temp2[temp[j+1]] = -1$ \\
                Set $D_g1_T.append(temp2)$ }}
		\caption{Generate $D_{g_1,T}$ for $g_1$}
	\label{algo:D_g1T}
\end{algorithm}
To satisfy the equation which conserve the energy, e.g. \eqref{disconserveEnergy} we will need the matrices
\begin{align*}
A_{g_1}^+ = c_p \rho A_i \mathcal{I}_r^+ \text{ and } A_{g_1}^- = c_p \rho A_i \mathcal{I}_r^-
\end{align*}
and set
\begin{align*}
A_{g_1,c1} &= \left( 0_{N \times \tilde{n}} | \eins_{N} \right)\\
A_{g_1,c3} &= \left( \eins_{N} | 0_{N \times 2N}  \right).
\end{align*}
Moreover we use $A_{g_1,c2} = 0_{N \times \tilde{n} + N}$ and set $A_{g_1,c2}[i, \tilde{n}^2_i] = 1$, where $\tilde{n}^2 = cumsum(n_i-2)$ and $i=1,...,N$. 
Setting the pressure at the end of the network and all given inlet temperatures is done by
\begin{align*}
D_{g1,AW} = \begin{pmatrix}
0_{(n_g1-(n_c+1)) \times 3N} \\
 \begin{pmatrix}0 & \hdots & 0 & 1	\end{pmatrix}_{3N}  \\
C^2 \quad 0_{n_c \times 2N}
\end{pmatrix} \text{ and } B_{AW}(t) = \begin{pmatrix}
0_{n_g1-(n_c+1)} \\
-p_N(t) \\
-T_out \\
\vdots
\end{pmatrix}.
\end{align*}
To simulate a network we need in addition
\begin{align*}
B_s(t) = \begin{pmatrix}
0_{n_p + n_{out}} \\
- T_{in}(t) \\
- p_{in}(t) \\
0_{n_{junc}+n_c +1}
\end{pmatrix} \text{ and } D_{g_1,AW1} = \begin{pmatrix}
 0_{(n_p+n_{out}) \times 3N}  \\
\begin{pmatrix}1 & 0 & \hdots & 0	\end{pmatrix}_{3N} \\
\begin{pmatrix}0_N & 1 & 0 & \hdots & 0	\end{pmatrix}_{3N} \\
0_{(n_{junc}+n_c+1) \times 3N}
\end{pmatrix}.
\end{align*}
Now we can summarize our matrices by
\begin{align*}
D_{g_{12}} &= D_{g1,p} + D_{g1,T} + D_{g1,AW} + D_{g1,AW1}, \\
B(t) &= B_{AW}(t) + B_{s}(t) 
\end{align*}
and set 
\begin{align*}
g_1(t,x(t),y(t)) = \begin{pmatrix}
D_{g_{12}} y + B(t) \\
A_{g1}^+((A_{g1,c1} x(t)).*(A_{g_1,c3}y(t))) + A_{g1}^-((A_{g1,c1} x(t)).*(A_{g1,c2} x(t)))
\end{pmatrix},
\end{align*}
where 
\begin{align*}
A_{g1,c1}, A_{g1,c2} \in \R^{2N \times \tilde{n}}, \, 
A_{g1,c3} \in \R^{2N \times 3N}, \, 
D_{g12} \in \R^{n_{g1} \times 3N}, \,
B(t) \in \R^{n_{g1}}.
\end{align*}
\paragraph{The function $g_2$}
\mbox{}\\
Now we define the matrices to get the function $g_2$ which depends only on $x$. We start with the equations for the conservation of mass and no mass is lost and set 
\begin{align*}
A_{g2,11} &= \begin{pmatrix}
0_{n_{junc} \times \tilde{n}}, & \mathcal{I}_r diag(A_1,...,A_N)
\end{pmatrix}, \\
A_{g2,12} &= \begin{pmatrix}
0_{n_c \times \tilde{n}}, & C^1 diag(\rho,...,\rho) - C^2 diag(\rho,...,\rho) 
\end{pmatrix}, \\
A_{g2,1} &= \begin{pmatrix}
A_{g2,11} \\
A_{g2,12} \\
0_{n_c \times (\tilde{n}+N)}
\end{pmatrix}
\end{align*}
which are already enough to satisfy \eqref{disconserveMass} and \eqref{disnoMassLost}. For the consumer demand equation \eqref{disconsumerDemand} we need
\begin{align*}
A_{g1,Q1} &= \begin{pmatrix}
0_{n_c \times \tilde{n}}, & C^1 diag(c_p\rho A_1,...,c_p\rho A_N) 
\end{pmatrix}, \\
A_{g2, 2} &= \begin{pmatrix}
0_{(n_{junc}+n_c) \times (\tilde{n}+N)} \\
A_{g1,Q1}
\end{pmatrix}, \\
A_{g2, 3} &= \begin{pmatrix}
0_{(n_{junc}+n_c) \times (\tilde{n}+N)} \\
A_{g1,Q2}
\end{pmatrix},\\
d_2 &= \begin{pmatrix}
0_{n_{junc} + n_c} \\
-T_{out_{n_c}}
\end{pmatrix},  \\
B_Q(t) &= \begin{pmatrix}
0_{n_{junc} + n_c} \\
-Q_1(t) \\
\vdots \\
-Q{_{n_c}}(t) 
\end{pmatrix},
\end{align*}
where $A_{g2,Q2}$ is calculated as described in Algorithm \ref{algo:g2_Q2}.
\begin{algorithm}
	\SetAlgoLined
    	\SetKwInOut{Output}{Return}
		\Output{matrix $A_{g2,Q2}$}
	Set $A1 = 0_{n_c \times \tilde{n}_1 -1}$ \\
	Set $A2 = \left( A1 \, | \, C^1[:,0]  \right)$
	\For{$i = 1,2, ..., N$}{
		Set $A1 = 0_{n_c \times \tilde{n}_i -1}$ \\ 
		Set $A2 = \left( A2, A1 \right)$ \\
		Set $A2 = \left( A2 \, | \, C^1[:,i]  \right)$ }
	Set $A_{g2,Q2} = \left( A2 \, | \, 0_{n_c \times N} \right)$
		\caption{Generate $A_{g2,Q2}$ for $g_2$}
	\label{algo:g2_Q2}
\end{algorithm}
All together we set
\begin{align*}
g_2(t,x(t)) = A_{g2,1} x(t) + A_{g_2,2} x(t) .* \left( A_{g2,3} x(t)  + d_2 \right) + B_{Q} (t),
\end{align*}
where
\begin{align*}
A_{g2,1}, A_{g_2,2}, A_{g_2,3} \in \R^{n_{junc}+2n_c \times \tilde{n}}, \, B_Q(t), d_2 \in R^{n_{junc}+2n_c}.
\end{align*}
\begin{remark}
\label{rem:singular_regular}
One can verify that $\frac{\partial g_1}{y}$ is regular and $\frac{\partial g_2}{y}$ is singular.
\end{remark}
%

\paragraph{Derivation of $f_1$, $f_2$, $g_1$ and $g_2$}
\mbox{}\\
Since we will need them we calculate all derivatives with respect to all arguments without t. Let us start with the function $f$ which is split into $f_1$ and $f_2$. As a reminder, here is the function
\begin{align*}
f_1(t, x, y) &= A_{f_1} x_2(t)  .* \left( A_{f_2}x_1(t) + D_{f_1}y(t) \right) +  A_{f_4} x_1(t) + d_1 
\end{align*}
and its derivatives are given by
\begin{align*}
\frac{d}{dx_1} f_1(t, x, y) &= (A_{f1}x_2) .* A_{f2} +  A_{f4}, \\
\frac{d}{dx_2} f_1(t, x, y) &= A_{f1} .*(A_{f2} x_1+D_{f1}y), \\
\frac{d}{dy} f_1(t, x, y) &= (A_{f1} x_2).*(D_{f1}).
\end{align*}
The function $f_2$ is given by
\begin{align*}
f_2(t, x(t), y(t)) = A_{f_3} x_2(t)  .* x_2(t)  +  D_{f_2} y(t) + d_2
\end{align*}
and its derivatives are given by
\begin{align*}
\frac{d}{dx_1} f_2(t, x, y) &= 0_{\tilde{n} \times N}, \\
\frac{d}{dx_2} f_2(t, x, y) &= A_{f3} .*x_2 + (A_{f3} x_2) .* 1_N,  \\
\frac{d}{dy} f_2(t, x, y) &= D_{f2}.
\end{align*}
Next we give the derivatives of $g_1$:
\begin{align*}
\frac{d}{dx} g_1(t,x,y,u) &= \begin{pmatrix}
A_{g11}.*(A_{g12}x+D_{g11}y)+(A_{g11} x).*(A_{g12}) \\
A_g1^+(A_{g1c1}.*(A_{g1c3} y)) +
A_g1^-(A_{g1c1}.*(A_{g1c2}x) + (A_{g1c1} x).*(A_{g1c2}))
\end{pmatrix}, \\
\frac{d}{dy} g_1(t,x,y,u) &= \begin{pmatrix}
D_{g12} + (A_{g11} x).*(D_{g11}) \\
A_{g1}^+((A_{g1c1} x).*(A_{g1c3}))
\end{pmatrix}, \\
\frac{d}{du} g_1(t,x,y,u) &= \begin{pmatrix}
B_C \\
0_{3*N-len(B_C) \times 3}
\end{pmatrix}.
\end{align*}
And finally the derivatives of $g_2$:
\begin{align*}
\frac{d}{dx} g_2(t, x) &=  A_{g21} + A_{g22}.*((A_{g23} x + d2)) +
               (A_{g22} x).*(A_{g23}), \\
\frac{d}{dy} g_2(t, x) &= 0_{(n_{junc}+2n_c) \times 3N}.\\
\end{align*}

\section*{Acknowledgments}
The authors acknowledge funding by the German Federal Ministry of Education and Research (BMBF) for the project DynOptHeat under the Förderkennzeichen (FKZ) 01LY1917B.

\bibliographystyle{plain}
\bibliography{references.bib}

\begin{thebibliography}{10}

\bibitem{Albaiz14}
Abdulaziz Albaiz.
\newblock High-order finite-difference discretization for steady-state
  convection-diffusion equation on arbitrary domain.
\newblock {\em Journal of Computational Physics}, 345(6):358--372, 2017.

\bibitem{Andersson2015}
Christian Andersson, Claus F\"uhrer, and Johan {\AA}kesson.
\newblock Assimulo: a unified framework for {ODE} solvers.
\newblock {\em Mathematics and Computers in Simulation}, 116(0):26--43, 2015.

\bibitem{Banda06}
Mapundi Banda, Michael Herty, and Axel Klar.
\newblock Coupling conditions for gas networks governed by the isothermal euler
  equations.
\newblock {\em NHM}, 1:295--314, 06 2006.

\bibitem{BangGutinGregory2008}
J{\o}rgen Bang-Jensen and Gregory~Z Gutin.
\newblock {\em Digraphs: theory, algorithms and applications}.
\newblock Springer Science \& Business Media, 2008.

\bibitem{bordin2016optimization}
Chiara Bordin, Angelo Gordini, and Daniele Vigo.
\newblock An optimization approach for district heating strategic network
  design.
\newblock {\em European Journal of Operational Research}, 252(1):296--307,
  2016.

\bibitem{Borsche19}
Raul Borsche, Matthias Eimer, and Norbert Siedow.
\newblock A local time stepping method for thermal energy transport in district
  heating networks.
\newblock {\em Applied Mathematics and Computation}, 353:215--229, 07 2019.

\bibitem{BuchvomSchropp}
Kathryn~Eleda Brenan, Stephen~L Campbell, and Linda~Ruth Petzold.
\newblock {\em Numerical solution of initial-value problems in
  differential-algebraic equations}.
\newblock SIAM, 1995.

\bibitem{Petzold}
Kathryn~Eleda Brenan, Stephen~L Campbell, and Linda~Ruth Petzold.
\newblock {\em Numerical solution of initial-value problems in
  differential-algebraic equations}.
\newblock SIAM, 1995.

\bibitem{burger2017survey}
Michael Burger and Matthias Gerdts.
\newblock A survey on numerical methods for the simulation of initial value
  problems with sdaes.
\newblock In {\em Surveys in Differential-Algebraic Equations IV}, pages
  221--300. Springer, 2017.

\bibitem{Colombo08}
Rinaldo~M. Colombo and Mauro Garavello.
\newblock On the cauchy problem for the p-system at a junction.
\newblock {\em SIAM Journal on Mathematical Analysis}, 39(5):1456--1471, 2008.

\bibitem{Diestel2017}
Reinhard Diestel.
\newblock {\em Graph Theory}.
\newblock Springer, 2017.

\bibitem{Domschke15}
Pia Domschke, Oliver Kolb, and Jens Lang.
\newblock Adjoint-based error control for the simulation and optimization of
  gas and water supply networks.
\newblock {\em Applied Mathematics and Computation}, 259:1003--1018, 05 2015.

\bibitem{Fuegenschuh15}
Armin F\"ugenschuh, Björn Gei{\ss}ler, Ralf Gollmer, Antonio Morsi, Jessica
  R\"ovekamp, Martin Schmidt, Klaus Spreckelsen, and Marc Steinbach.
\newblock {\em Chapter 2: Physical and technical fundamentals of gas networks},
  pages 17--43.
\newblock 03 2015.

\bibitem{gerdts2011optimal}
Matthias Gerdts.
\newblock {\em Optimal control of ODEs and DAEs}.
\newblock Walter de Gruyter, 2011.

\bibitem{Herty2006}
Michael Herty.
\newblock Modeling, simulation and optimization of gas networks with
  compressors.
\newblock {\em Networks and Heterogeneous Media}, 2(1):81--97, 2006.

\bibitem{Sundials}
Alan~C Hindmarsh, Peter~N Brown, Keith~E Grant, Steven~L Lee, Radu Serban,
  Dan~E Shumaker, and Carol~S Woodward.
\newblock {SUNDIALS}: Suite of nonlinear and differential/algebraic equation
  solvers.
\newblock {\em ACM Transactions on Mathematical Software (TOMS)},
  31(3):363--396, 2005.

\bibitem{Koecher00}
Ralf K\"ocher.
\newblock {B}eitrag zur {B}erechnung und {A}uslegung von {F}ernw\"armenetzen.
\newblock 01 2000.

\bibitem{Krug19}
Richard Krug, Volker Mehrmann, and Martin Schmidt.
\newblock Nonlinear optimization of district heating networks.
\newblock {\em Optimization and Engineering}, 22(1):1--37, 2021.

\bibitem{kunkelBook}
Peter Kunkel and Volker Mehrmann.
\newblock {\em Differential-algebraic equations: analysis and numerical
  solution}, volume~2.
\newblock European Mathematical Society, 2006.

\bibitem{odes}
Benny Malengier, Pavol Ki{\v{s}}on, James Tocknell, Claas Abert, Florian
  Bruckner, and Marc-Antonio Bisotti.
\newblock {ODES}: a high level interface to {ODE} and {DAE} solvers.
\newblock {\em The Journal of Open Source Software}, 3(22):165, Feb 2018.

\bibitem{mehrmann2012index}
Volker Mehrmann.
\newblock Index concepts for differential-algebraic equations.
\newblock {\em Encyclopedia of Applied and Computational Mathematics},
  1:676--681, 2012.

\bibitem{Qiu18}
Yue Qiu, Sara Grundel, Martin Stoll, and Peter Benner.
\newblock Efficient numerical methods for gas network modeling and simulation.
\newblock 07 2018.

\bibitem{rein2019model}
Markus Rein, Jan Mohring, Tobias Damm, and Axel Klar.
\newblock Model order reduction of hyperbolic systems at the example of
  district heating networks, 2021.

\bibitem{schwarz2018new}
Diana~Est{\'e}vez Schwarz and Ren{\'e} Lamour.
\newblock A new approach for computing consistent initial values and taylor
  coefficients for daes using projector-based constrained optimization.
\newblock {\em Numerical Algorithms}, 78:355--377, 2018.

\bibitem{Schwarz2000Consistent}
Diana~Estévez Schwarz.
\newblock {\em Consistent initialization for index-2 differential algebraic
  equations and its application to circuit simulation}.
\newblock PhD thesis, Humboldt-Universität zu Berlin,
  Mathematisch-Naturwissenschaftliche Fakultät II, 2000.

\bibitem{Hairer2}
Gerhard Wanner and Ernst Hairer.
\newblock {\em Solving ordinary differential equations II}.
\newblock Springer Berlin Heidelberg, 1996.

\end{thebibliography}

\printindex

\end{document}